\documentclass[12pt, a4paper]{article}

\usepackage[T1]{fontenc}
\usepackage[latin1]{inputenc}
\usepackage[english]{babel}
\usepackage{amsmath}
\usepackage{textcomp}
\usepackage{graphicx}
\usepackage{amsfonts}
\newtheorem{theorem}{Theorem}[section]
\newtheorem{remark}{Remark}[section]
\newtheorem{lemma}{Lemma}[section]
\newtheorem{corollary}{Corollary}[section]
\makeatletter

\@addtoreset{equation}{section}
\makeatother

\title{KPP reaction-diffusion systems with loss inside a cylinder: convergence toward the problem with Robin boundary conditions}
\author{Thomas Giletti \thanks{Université Aix-Marseille III, LATP, Faculté des Sciences et Techniques, Avenue Escadrille Normandie-Niemen, F-13397 Marseille Cedex 20, France; thomas.giletti@etu.univ-cezanne.fr}}

\begin{document}

\maketitle

\begin{abstract}
We consider in this paper a reaction-diffusion system under a KPP hypothesis in a cylindrical domain in the presence of a shear flow. Such systems arise in predator-prey models as well as in combustion models with heat losses. Similarly to the single equation case, the existence of a minimal speed $c^*$ and of traveling front solutions for every speed $c>c^*$ has been shown both in the cases of heat losses distributed inside the domain or on the boundary. Here, we deal with the accordance between the two models by choosing heat losses inside the domain which tend to a Dirac mass located on the boundary. First, using the characterizations of the corresponding minimal speeds, we will see that they converge to the minimal speed of the limiting problem. Then, we will take interest in the convergence of the traveling front solutions of our reaction-diffusion systems. We will show the convergence under some assumptions on those solutions, which in particular can be satisfied in dimension 2.
\end{abstract}

\section{Introduction and main results}\label{sec:intro}

\subsection*{The models and their background}

We consider reaction-diffusion-advection systems in a cylindrical domain $\Omega\! =\! \mathbb{R}_x\! \times \! \omega_y \! \subset \mathbb{R}^d$ where $\omega$ is a smooth bounded domain of $\mathbb{R}^{d-1}$. The existence and qualitative properties of the solutions of such problems have been extensively studied over the years both in the single-equation \cite{berestycki5,berestycki6,kiselev1} and the two-equations \cite{hamel-quenching,gordon-periodic,hamel-nonadiabatic,hamel-adiabatic,roques2} cases. Those references describe various situations, in dimension 1 or more, within homogeneous or heterogeneous framework, and with various assumptions on the nonlinearities arising in the system, such as KPP or ignition hypotheses. For large reviews of this mathematical area, we refer the reader to \cite{berestycki-hamel,berestycki4,xin3}. This variety reflects the diversity of the processes which can lead to such systems, ranging from chemical and biological to combustion contexts \cite{murray1,williams}.

In particular, in this paper, the issue at stake is the accordance between two reaction-diffusion-advection problems in $\Omega$ with heat losses. For a better knowledge of those models, we will refer the reader to \cite{hamel-quenching,giletti1,hamel-nonadiabatic}, where they have been introduced.  Note that here, we chose to invoke the "combustion" terminology, hence the term "heat loss". We will also refer to the unknowns functions $Y$ and $T$ as, respectively, the combustant concentration and the temperature.

Let us now present the two models. First, when the heat loss (denoted by $h$) can take place in the whole domain $\Omega$, we consider the following system, described in \cite{giletti1}:
\begin{equation}\label{eqn:sysh}
\left\{
\begin{array}{l}
T_t + u(y) T_x= \Delta T + f(y,T)Y - h(y,T),\vspace{3pt} \\
Y_t + u(y) Y_x= \mbox{Le}^{-1} \Delta Y - f(y,T)Y,\\
\end{array}
\right.
\end{equation}
with Neumann boundary conditions
\begin{equation}\label{eqn:neumann}
\frac{\partial T}{\partial n} = \frac{\partial Y}{\partial n} = 0 \mbox{ on } \partial \Omega,
\end{equation}
where $n$ denotes the outward unit normal on $\partial \Omega$.

Then, in the case of a heat loss (denoted by $qT$) on the boundary, we consider the following system, described in \cite{hamel-quenching,hamel-nonadiabatic}:
\begin{equation}\label{eqn:sysrobin}
\left\{
\begin{array}{l}
T_t + u(y) T_x= \Delta T + f(y,T)Y ,\vspace{3pt} \\
Y_t + u(y) Y_x= \mbox{Le}^{-1} \Delta Y - f(y,T)Y,\\
\end{array}
\right.
\end{equation}
with Robin boundary conditions
\begin{equation}\label{eqn:robin}
\left\{
\displaystyle
\begin{array}{l}
\displaystyle \vspace{1mm}
\frac{\partial T}{\partial n} + q T = 0 \mbox{ on } \partial \Omega,\vspace{3pt}  \\
\displaystyle \vspace{1mm}
\frac{\partial Y}{\partial n} = 0 \mbox{ on } \partial \Omega. \\
\end{array}
\right.
\end{equation}

Note that in both systems (\ref{eqn:sysh}) and (\ref{eqn:sysrobin}), the Lewis number $\mbox{Le}$ is the ratio of the thermal diffusivity and the diffusivity of the reactant, and may be an arbitrary positive number. Besides, $f(.,T)Y$ is the reaction term, raising the temperature while consuming the combustant, and $u$ is the shear flow of the medium, which will be assumed to have zero average and does not depend on the $x$-variable: that is, the flow is invariant along the cylinder, and divergence free.

As mentioned before, those systems have a wide range of applications, and also describe predator-prey situations \cite{murray1}. The unknowns $T$ and $Y$ would then be replaced by the density of two species, $U$ and $V$, where the former is the predator and the latter the prey. In this context, our interest in traveling waves is a way to study the invasion of the prey-populated medium by the predator. The "heat loss" would then be interpreted as the death rate or as a saturation effect for the species $U$. In the former, the death rate of the predator $U$ may be caused by a human influence, which may intervene inside the medium or on its boundary. In the latter, the saturation may be caused by an intra-species competition. In both cases, it can be reasonable to assume that the intrinsic death rate and the saturation for the prey $V$ are negligible compared to the other parameters.

As we are interested in traveling front solutions of those problems, we look for solutions of the form $T(t,x,y)=\tilde{T} (x-ct,y)$, $Y(t,x,y)=\tilde{Y} (x-ct,y)$ for some $c \in \mathbb{R}$, and such that
\begin{equation}\label{eqn:condinfty}
\left\{
\begin{array}{l}
\tilde{T} (+\infty ,.)=0, \ \tilde{Y}(+\infty ,.)=1 , \vspace{3pt} \\
\tilde{T}_x (-\infty ,.)=\tilde{Y}_x (-\infty ,.)=0 .\\
\end{array}
\right.
\end{equation}
Physically, that means that we search solutions such that the right side of the front is a cold region with a strong presence of the combustant (or, in a biological context, is populated only by the prey), while the left side of the front is left free. Moreover, to make physical sense, we impose on the solutions to verify $T>0$ and $0<Y<1$ (the inequalities are strict in order to avoid trivial solutions).
\subsection*{Assumptions, notations and known results}

Before we enounce the main results of \cite{giletti1} and \cite{hamel-quenching,hamel-nonadiabatic}, we remind the main hypotheses. We assume first that $u \in C^{0 ,\alpha} (\overline{\omega})$ (for some $\alpha >0$) and that, as said before, it has zero average:
$$\int_\omega u(y) dy =0.$$
Moreover, the heat loss coefficient $q$ is assumed, as in \cite{hamel-quenching,hamel-nonadiabatic}, to be a positive constant, although the results in the mentioned papers and in this article could be generalized to smooth positive functions on $\partial \omega$.

The functions $f$ and $h$ are in $C^1(\overline{\omega}\times[0,+\infty);\mathbb{R})$ and there exists $s_0 >0$ such that the sets of functions $(f(y,.))_{y\in \overline{\omega}}$ and $(h(y,.))_{y\in \overline{\omega}}$ are bounded in $C^{1,\alpha} ( [0,s_0) ; \mathbb{R}) $. Lastly, $f$ verifies:
\begin{equation}\label{condf}
\begin{array}{l}
\displaystyle
\displaystyle f(.,0)=0 <  f(.,T) \leq \frac{\partial f}{\partial T}(.,0)T , \ \frac{\partial f}{\partial T} \geq 0 \mbox{ for all } T>0, \ \mbox{and } f(., +\infty )=+\infty ;\\
\end{array}
\end{equation}
and $h$ satisfies:
\begin{equation}\label{condh}
\left\{
\begin{array}{c}
\displaystyle
h(.,0)=0 \leq \frac{\partial h}{\partial T}(.,0)T \leq h(.,T) \leq KT \mbox{ for all } T \geq 0 \mbox{ and some } K>0, \vspace{3pt} \\
\displaystyle \int_{\omega} \frac{\partial h}{\partial T}(y,0) dy > 0.
\end{array}
\right.
\end{equation}
The KPP-type hypotheses will allow us to determine the behavior of systems (\ref{eqn:sysh})-(\ref{eqn:neumann}) and (\ref{eqn:sysrobin})-(\ref{eqn:robin}) by comparisons with the linearized problems. The positivity of the integral of $\frac{\partial h}{\partial T} (.,0)$ insure that the heat-loss is non trivial. For a more precise information about the role of those hypotheses, we refer to  \cite{hamel-quenching,giletti1,hamel-nonadiabatic}.

Next, in order to characterize the minimal speeds of the solutions of the two systems (\ref{eqn:sysh})-(\ref{eqn:neumann}) and (\ref{eqn:sysrobin})-(\ref{eqn:robin}), we introduce some eigenvalue problems that arise from the linearized systems ahead of the front (that is, with $T=0$ and $Y=1$). For $\lambda \in \mathbb{R}$, let~$(\mu_h (\lambda), \phi_{h,\lambda} )$ be the principal eigenvalues and eigenfunctions of the following system:
\begin{equation}\label{eqn:principaleigenvalueh}
\left\{
\begin{array}{rcll}
\displaystyle -\Delta \phi_{h ,\lambda} - \lambda u(y) \phi_{h,\lambda} + \left(\frac{\partial h}{\partial T}  (y,0) -\frac{\partial f}{\partial T} (y,0)\right)  \phi_{h,\lambda} & = & \mu_{h} (\lambda ) \phi_{h,\lambda} & \mbox{ in } \omega ,\vspace{3pt} \\
\displaystyle \frac{\partial \phi_{h,\lambda} }{\partial n} & = & 0 & \mbox{ on } \partial \omega .\\
\end{array}
\right.
\end{equation}
That is, $\mu_h (\lambda )$ is the unique eigenvalue of (\ref{eqn:principaleigenvalueh}) that corresponds to a positive eigenfunction~$\phi_{h,\lambda}$. Nonnegative solutions of the form $T(t,x,y)=\phi (y) e^{-\lambda (x-ct)}$ of the linearized system (\ref{eqn:sysh})-(\ref{eqn:neumann}) with $T=0$ and $Y=1$ only exist if $\phi=\phi_\lambda$ and $\mu_{h,f} (\lambda)=\lambda^2 -c\lambda$. For the conditions $(\ref{eqn:condinfty})$ ahead of the front to be satisfied, the real number $\lambda$ must be positive, thus the need of assumptions on $\mu_{h,f}$ and $c$, which will be enounced later in this paper and will in fact guarantee the existence of traveling waves (see Theorem~\ref{th:thremind} below). Similarly, let also $(\nu_q  (\lambda), \psi_{q,\lambda} )$ be the principal eigenvalues and eigenfunctions of the system:
\begin{equation}\label{eqn:principaleigenvaluerobin}
\left\{
\begin{array}{rcll}
\displaystyle -\Delta \psi_{q,\lambda} - \lambda u(y) \psi_{q,\lambda} -\frac{\partial f}{\partial T} (y,0)  \psi_{q,\lambda} & = & \nu_q  (\lambda) \psi_{q,\lambda} & \mbox{ in } \omega ,\vspace{3pt} \\
\displaystyle \frac{\partial \psi_{q,\lambda} }{\partial n} + q \psi_{q,\lambda} & = & 0 & \mbox{ on } \partial \omega .\\
\end{array}
\right.
\end{equation}
This system arises from the linearization of (\ref{eqn:sysrobin})-(\ref{eqn:robin}) with $T=0$ and $Y=1$.
We can normalize $\phi_{h,\lambda}$ and $\psi_{q,\lambda}$ in $L^2 (\omega)$ norm, that is:
\begin{equation}\label{eqn:L2norm}
\| \phi_{h ,\lambda} \|_2 =
\| \psi_{q,\lambda} \|_2 =1 \ .
\end{equation}
We now enounce some properties of $\nu_q  (\lambda)$ and $\mu_h (\lambda )$ that will be needed throughout this paper. First, we remind that under the $L^2$-normalization (\ref{eqn:L2norm}):
\begin{equation}\label{eqn:principaleigenvalueinfmu}
\begin{array}{rcll}
\displaystyle
\hspace{15mm}\mu_{h} (\lambda ) & = & \displaystyle  \inf_{\phi \in H^1 (\omega) , \| \phi \|_2 = 1} \displaystyle \left[ \int_\omega \displaystyle | \nabla \phi (y) |^2 dy - \lambda \int_\omega u(y) \phi ^2 (y) dy \right.
 \vspace{3pt}\\
 & & \displaystyle \left. \hspace{3cm} \displaystyle + \int_\omega \displaystyle \left( \displaystyle \frac{\partial h}{\partial T}  (y,0) -\frac{\partial f}{\partial T} (y,0)\right) \phi ^2 (y) dy \right] , \vspace{3pt}\\
\displaystyle
& = & \displaystyle \int_\omega \displaystyle | \nabla \phi_{h,\lambda} (y) |^2 dy - \lambda \int_\omega u(y) \phi_{h,\lambda}^2 (y) dy \vspace{3pt} \\
& & \hspace{3cm} \displaystyle + \int_\omega \left( \frac{\partial h}{\partial T}  (y,0) -\frac{\partial f}{\partial T} (y,0)\right) \phi_{h,\lambda} ^2 (y) dy ,
\end{array}
\end{equation}
\begin{equation}\label{eqn:principaleigenvalueinfnu}\displaystyle
\begin{array}{rcl}
\displaystyle
\nu_q (\lambda ) & = & \displaystyle \inf_{\phi \in H^1 (\omega) , \| \phi \|_2 = 1} \left[ \displaystyle \int_\omega | \nabla \phi (y) |^2 dy - \lambda \int_\omega u(y) \phi^2 (y) dy  \vspace{3pt} \right.\\
& & \hspace{3cm} \left. \displaystyle + \int_{\partial \omega} q \phi^2 - \int_\omega \frac{\partial f}{\partial T} (y,0) \phi^2 (y) dy \right] , \vspace{3pt}  \\
\displaystyle & = & \displaystyle \int_\omega | \nabla \psi_{q,\lambda} (y) |^2 dy - \lambda \int_\omega u(y) \psi_{q,\lambda} ^2 (y) dy \vspace{3pt}\\
& & \hspace{3cm} \displaystyle + \int_{\partial \omega} q \psi_{q,\lambda}^2 - \int_\omega \frac{\partial f}{\partial T} (y,0) \psi_{q,\lambda}^2 (y) dy  .
\end{array}
\end{equation}
It follows from the above that the functions $\mu_h (\lambda )$ and $\nu_q (\lambda )$ are concave, as an infimum of affine functions of $\lambda$. Lastly, elementary calculations lead to, for all $\lambda \in \mathbb{R}$:
\begin{equation}\label{eqn:derivative}
\left\{
\begin{array}{rcl}
\displaystyle \mu '_h (\lambda )  & = & - \displaystyle \int_\omega \displaystyle u(y) \phi_{h,\lambda}^2 (y) dy, \vspace{3pt} \\
\displaystyle \nu_q ' (\lambda )  & = &  - \displaystyle \int_\omega \displaystyle u(y) \psi_{q,\lambda}^2 (y) dy.
\end{array}
\right.
\end{equation}
We can now express the minimal speeds of problems (\ref{eqn:sysh})-(\ref{eqn:neumann}) and (\ref{eqn:sysrobin})-(\ref{eqn:robin}) with conditions at infinity (\ref{eqn:condinfty}) as:
\begin{equation}\label{eqn:minspeed}
\left\{
\begin{array}{rcl}
\displaystyle c_h^* & =  & \displaystyle  \min \{ c \in \mathbb{R} / \exists \lambda >0 , \mu_h (\lambda) = \lambda^2 - c\lambda \}, \vspace{3pt}\\
\displaystyle c_q^* & = & \displaystyle \min \{ c \in \mathbb{R} / \exists \lambda >0 , \nu_q (\lambda ) = \lambda^2 - c\lambda \}.
\end{array}
\right.
\end{equation}
We assume in this paper that:
\begin{equation*}
\left\{
\begin{array}{rcl}
\mu_h (0)& <& 0, \vspace{3pt} \\
\nu_q (0)& < & 0.
\end{array}
\right.
\end{equation*}
Then, by concavity of $\mu_h$ and $\nu_q$ with respect to $\lambda$, the minima in (\ref{eqn:minspeed}) are well defined.

The known existence results for problems $(\ref{eqn:sysh})$-$(\ref{eqn:neumann})$ and $(\ref{eqn:sysrobin})$-$(\ref{eqn:robin})$ are summed up in the following theorem:
\begin{theorem}[\cite{hamel-quenching,giletti1,hamel-nonadiabatic}]\label{th:thremind} Under the above hypotheses, we have that:

$(a)$ If $\mu_h (0) <0$, there exists a traveling front solution of $(\ref{eqn:sysh})$-$(\ref{eqn:neumann})$ and $(\ref{eqn:condinfty})$ with speed~$c$ if and only if $c \geq c_h^*$ and $c>0$,

$(b)$ If $\nu_q (0) <0$ and $c > \max ( c_q^* ,0)$, then there exists a traveling front solution of $(\ref{eqn:sysrobin})$-$(\ref{eqn:robin})$ and $(\ref{eqn:condinfty})$ with speed $c$. Conversely, if there exists a traveling front solution of $(\ref{eqn:sysrobin})$-$(\ref{eqn:robin})$ and $(\ref{eqn:condinfty})$ with speed $c$, then $c \geq c_q^*$ and $c>0$.
\end{theorem}
\subsection*{Main results}
In this paper, we first prove the following result on the convergence of the minimal speeds:
\begin{theorem}\label{th:main} Let $(h_k)_{k \in \mathbb{N}}$ be a sequence of functions verifying $(\ref{condh})$ with $h=h_k$ for all~$k \in \mathbb{N}$, and such that:
\begin{equation}\label{eqn:hypbord}
\exists \varepsilon_k \searrow 0 \mbox{ such that } \frac{\partial h_k}{\partial T}  (.,0) \rightarrow 0 \mbox{ uniformly in }\omega \backslash (\partial \omega + \overline{B} (0,\varepsilon_k)),
\end{equation}
\begin{equation}\label{eqn:hypO}
\displaystyle \varepsilon_k \left\| \frac{\partial h _k}{\partial T}  (.,0) \right\|_{L^\infty (\omega )} = O (1),
\end{equation}
\begin{equation}\label{eqn:hyp}
g(\sigma) := \int_0^1 \varepsilon_k \frac{\partial h _k}{\partial T}  (\sigma -\varepsilon_k s n(\sigma),0) ds \ \longrightarrow \ q \ \mbox{ uniformly in } \sigma \in \partial \omega .
\end{equation}

Then $\mu_{h_k} \rightarrow \nu_{q}$ locally uniformly and for any $\lambda \in \mathbb{R}$, the sequence $(\phi_{h_k ,\lambda})_{k \in \mathbb{N}}$ of the $L^2$-normalized principal eigenfunctions of $(\ref{eqn:principaleigenvalueh})$ is bounded in $H^1 (\omega )$ and converges to the principal eigenfunction $\psi_{q , \lambda}$ of $(\ref{eqn:principaleigenvaluerobin})$ strongly in $L^2 (\omega)$ and weakly in $H^1 (\omega)$.

Furthermore, if $\nu_q (0)<0$, then $c^*_{h_k} \rightarrow c_q^*$.
\end{theorem}

Theorem \ref{th:main} will rely on the Lemma~\ref{lemma1}, that is the uniform convergence of the sequence $(\frac{\partial h_k}{\partial T} (.,0))_{k \in \mathbb{N}}$ toward the Dirac mass $q \delta_{\partial \omega}$ on any sequence of functions bounded in~$H^1 (\omega)$. This is in fact the important assumption on the sequence $(h_k)_{k \in \mathbb{N}}$, but we chose in Theorem~\ref{th:main} a more convenient hypothesis to give a better view of the admissible $h_k$. For instance, the $h_k$ defined as $$\displaystyle h_k (y,T) = \frac{\partial h_k}{\partial T} (y,0) \ T = 3 k^3 \left(\min\left(0, d(y,\partial \omega) - \frac{1}{k}\right)\right)^2 \ T$$
verify the correct assumptions with $\varepsilon_k = \frac{1}{k}$ and $q=1$. Since we assumed that $h_k$ belongs to $C^1~(\overline{\omega}~\times~[0~,+\infty~ );\mathbb{R})$, we can not define $\frac{\partial h_k}{\partial T} (.,0)$ as a constant or an affine function of~$d(.,\partial \omega)$ near $\partial \omega$, and $0$ elsewhere. Still, we can also define $\frac{\partial h_k}{\partial T} (.,0)$ as a constant or an affine function of $d(.,\partial \omega)$ near $\partial \omega$ and extend it so that $\frac{\partial h_k}{\partial T} (.,0) \in C^1 (\overline{\omega})$ and $\frac{\partial h_k}{\partial T} (.,0) =0$ outside of a neighborhood of $\partial \omega$. More generally, we can easily use standard approximations of a Dirac mass in dimension 1 to generate suitable sequences for Theorem~\ref{th:main}.

We also want to study the convergence of the traveling front solutions of the problem (\ref{eqn:sysh})-(\ref{eqn:neumann}) with the conditions at infinity (\ref{eqn:condinfty}), when $h$ converges to a Dirac mass with the same assumptions as in Theorem~\ref{th:main}. We recall that $(T,Y)=(T(x,y),Y(x,y))$ is a traveling front solution with speed $c \in \mathbb{R}$ of (\ref{eqn:sysh}) when:
\begin{equation}\label{eqn:syshfront}
\left\{
\begin{array}{l}
\Delta T + (c-u(y))T_x+ f(y,T)Y - h(y,T) =0 ,\vspace{3pt}\\
\mbox{Le}^{-1} \Delta Y + (c-u(y)) Y_x - f(y,T)Y =0.\\
\end{array}
\right.
\end{equation}
with Neumann boundary conditions (\ref{eqn:neumann}) and conditions at infinity (\ref{eqn:condinfty}). Similarly, $(T,Y)$ is a traveling front solution with speed $c \in \mathbb{R}$ of (\ref{eqn:sysrobin}) when:
\begin{equation}\label{eqn:sysrobinfront}
\left\{
\begin{array}{l}
\Delta T + (c-u(y))T_x+ f(y,T)Y =0 ,\vspace{3pt}\\
\mbox{Le}^{-1} \Delta Y + (c-u(y)) Y_x - f(y,T)Y =0.\\
\end{array}
\right.
\end{equation}
with Robin boundary conditions (\ref{eqn:robin}) and conditions at infinity (\ref{eqn:condinfty}).
To ensure the accordance between our two models of the heat-loss, we want to show the convergence of the solutions of (\ref{eqn:syshfront})-(\ref{eqn:neumann}) and (\ref{eqn:condinfty}) for some speed $c$ to the solutions of (\ref{eqn:sysrobinfront})-(\ref{eqn:robin}) and (\ref{eqn:condinfty}) with the same speed, when $h$ is replaced by a sequence $(h_k)_{k \in \mathbb{N}}$ converging to the Dirac mass $q \delta_{\partial \omega}$ with the same assumptions as in Theorem \ref{th:main}. The main difficulty is the lack of bounds on $h_k$, and thus the lack of estimates on the corresponding sequence of temperatures $T_k$. In particular, we would need $H^1$ estimates on the temperatures in order to use Lemma \ref{lemma1} on the convergence of the $h_k$ on any $H^1$~-bounded sequence toward a Dirac mass, that is the lemma we use for the convergence of the eigenvalue problems. Another difficulty is the impossibility to use Harnack inequality, so that even when our sequence of solutions converges, the limit may be trivial: for instance, the temperatures might tend to be concentrated on a single point, hence a regular limit could only be zero.

Therefore, the general case is still open. Here, we will only consider particular solutions of (\ref{eqn:syshfront})-(\ref{eqn:neumann}) and (\ref{eqn:condinfty}), which satisfy some exponential bounds from below and above, in order to overcome the above difficulties. Then, the bounds are proved to hold in dimension~$d=2$, leading to the desired convergence result in this case.

\begin{theorem}\label{th:main2}
Let $(h_k)_{k \in \mathbb{N}}$ be a sequence of functions verifying $(\ref{eqn:hypbord})$, $(\ref{eqn:hypO})$, $(\ref{eqn:hyp})$ and~$(\ref{condh})$ with $h=h_k$ for all $k \in \mathbb{N}$, and such that $h_k (y,.)$ is linear for all $y \in \overline{\omega}$ and $k \in \mathbb{N}$. Let also $(T_k , Y_k)$ a sequence of non-trival solutions of problem $(\ref{eqn:syshfront})$-$(\ref{eqn:neumann})$ with $h=h_k$, $c>c_q^*$, and verifying the conditions at infinity $(\ref{eqn:condinfty})$, and let $(\lambda_k)_{k \in \mathbb{N}}$ be a sequence of positive real numbers such that for any $k \in \mathbb{N}$, we have: $$\lambda_k ^2 - c \lambda_k = \mu_{h_k} (\lambda_k ).$$
We assume that there exist $0 < \Lambda_1 < \Lambda_2$ and $C_1 , C_2 , C_3 >0$ such that for all $k \in \mathbb{N}$ and~$(x,y) \in \overline{\Omega}$:
\begin{eqnarray}\label{eqn:Tnabove}
T_k (x,y) < C_1 e^{-\lambda_k x} ,
\end{eqnarray}
\begin{eqnarray}\label{eqn:Tnbelow}
\max (0, C_2 e^{-\Lambda_1 x} - C_3 e^{-\Lambda_2 x}) < T_k (x,y),
\end{eqnarray}
Then up to extraction of a subsequence, $(T_k , Y_k)$ converges weakly in $H^1_{loc} (\overline{\Omega} )$ and strongly in $L^2_{loc} (\overline{\Omega} )$ to a non trivial solution $(T,Y)$ of problem $(\ref{eqn:sysrobinfront})$-$(\ref{eqn:robin})$ and $(\ref{eqn:condinfty})$.
\end{theorem}

\begin{remark}In this theorem, we added the assumption that $h_k$ is linear in the $T$-variable. Indeed, in Theorem \ref{th:main} where we only considered the eigenvalue problem $(\ref{eqn:principaleigenvalueh})$, we only made assumptions on $\frac{\partial h_k}{\partial T} (.,0)$. Here, we need to make sure that the term "$h_k (y,T)$" in our equation $(\ref{eqn:syshfront})$ will tend to $q \delta_{\partial \omega} T$, hence the linearity assumption.
\end{remark}
The hypotheses of Theorem \ref{th:main2} may in fact be weakened. For instance, the hypothesis~(\ref{eqn:Tnbelow}) could be replaced by any function positive on a non trivial set in $\Omega$. The choice of those exponential bounds in fact come from the sub and super-solutions that have been used in \cite{giletti1} to construct solutions of (\ref{eqn:syshfront})-(\ref{eqn:neumann}) and (\ref{eqn:condinfty}). More precisely, we know that there exist, for each $k$ and $c >c_{h_k}^*$, solutions of (\ref{eqn:syshfront})-(\ref{eqn:neumann}) and (\ref{eqn:condinfty}) with $h=h_k$ that fulfill similar exponential bounds. The issue is then to make those bounds independent of $k \in \mathbb{N}$, in order to exhibit some particular solutions satisfying the assumptions of Theorem~\ref{th:main2}. In fact, the construction of sub and super-solutions aforementioned and used in~\cite{giletti1}, rely on some strong estimates on the principal eigenfunctions of (\ref{eqn:principaleigenvalueh}), which can be made independently of~$k \in \mathbb{N}$ only in dimension 2 ($d=2$), where $H^1 (\omega )$-estimates imply $C^{0, 1/2} (\omega)$-estimates.

The discussion above will lead to the following corollary of Theorem \ref{th:main2}:

\begin{corollary}\label{th:main22}
Let $(h_k)_{k \in \mathbb{N}}$ be a sequence of functions verifying $(\ref{eqn:hypbord})$, $(\ref{eqn:hypO})$, $(\ref{eqn:hyp})$ and~$(\ref{condh})$ with $h=h_k$ for all $k \in \mathbb{N}$ and such that $h_k (y,.)$ is linear for all $y \in \overline{\omega}$ and~$k \in \mathbb{N}$. In dimension $d=2$, up to extraction of some subsequence, there exists a sequence of solutions of problem $(\ref{eqn:syshfront})$-$(\ref{eqn:neumann})$ and $(\ref{eqn:condinfty})$ with $h=h_k$, $c> \max (0,c_q^*)$, that converges weakly in $H^1_{loc} (\overline{\Omega})$ and strongly in $L^2_{loc} (\overline{\Omega})$ to a non trivial solution $(T,Y)$ of problem $(\ref{eqn:sysrobinfront})$-$(\ref{eqn:robin})$ and~$(\ref{eqn:condinfty})$.
\end{corollary}
In spite of all the difficulties aforementioned, we think that this result is only a first step and in fact holds in a more general case. In particular, we hope that the study of the exponential behavior of any solution of (\ref{eqn:syshfront})-(\ref{eqn:neumann}) and (\ref{eqn:condinfty}), which will be a subject of interest in a forthcoming paper, will allow us to apply Theorem~\ref{th:main2} to a larger set of solutions.

\subsubsection*{Plan of the paper}

Theorem \ref{th:main} will be proved in Section \ref{sec:cveigen}. The use of Lemma \ref{lemma1} on the principal eigenfunctions of problem (\ref{eqn:principaleigenvalueh}) and (\ref{eqn:principaleigenvaluerobin}) will allow us to prove the locally uniform convergence of the eigenvalues $(\mu_k (\lambda ))_{k \in \mathbb{N}}$ toward~$\nu_q (\lambda )$. Lastly, we will end the proof by showing the convergence of the minimal speeds.

Theorem \ref{th:main2} will be proved in Section \ref{sec:cvsol}. We will first show Lemma \ref{lemmaY} which gives a uniform exponential bound from below on the sequence $(Y_k)_k$ near $+\infty$. Then, the bounds from above will allow us to obtain $H^1_{loc} (\overline{\Omega})$ estimates on the sequence $(T_k ,Y_k)_{k \in \mathbb{N}}$ and thus its convergence toward a pair $(T,Y)$. Then, the same lemma as in the proof of Theorem~\ref{th:main} will imply that $(T,Y)$ is a solution of (\ref{eqn:sysrobinfront})-(\ref{eqn:robin}). The fact that it is non trivial will immediately follow from (\ref{eqn:Tnbelow}) and Lemma~\ref{lemmaY}, and so will the behavior of $(T,Y)$ near $+\infty$. Lastly, the behavior of $(T,Y)$ on the left, near $-\infty$, will be proved using~(\ref{eqn:Tnabove}) and a lemma from~ \cite{hamel-nonadiabatic}, stating the boundedness of a solution of (\ref{eqn:sysrobinfront})-(\ref{eqn:robin}) when it is bounded from above by an exponential of the form $e^{-\lambda x}$ with $\lambda^2 - c\lambda = \nu_q (\lambda )$ (we include its proof at the end of Section \ref{sec:cvsol} for the sake of completeness).

The Section \ref{sec:cvsol2} will deal with the proof of Corollary \ref{th:main22}, although for convenience, we will refer the reader to \cite{giletti1} for the precise proof of the existence of solutions between the introduced sub and super-solutions.

\section{Convergence of the principal eigenvalue problems and of the minimal speeds}\label{sec:cveigen}

We deal in this section with the proof of Theorem \ref{th:main}. As mentioned before, we begin by a lemma before we study the convergence of the principal eigenvalues and minimal speeds.

\subsection{A useful lemma}

We prove here a lemma, that is the convergence of the sequence $(\frac{\partial h_k}{\partial T} (.,0))_{k \in \mathbb{N}}$ toward the Dirac mass $q \delta_{\partial \omega}$ in the following sense:

\begin{lemma}\label{lemma1}
Let $(\phi_k)_{k \in \mathbb{N}}$ be a bounded sequence of functions in $H^1 (\omega)$.
Then, up to extraction of some subsequence, the sequence converges weakly in $H^1 (\omega)$ and strongly in $L^2 (\omega )$ to a function $\phi$ such that:
$$\lim_n \int_\omega \frac{\partial h_k}{\partial T}  (.,0)  \phi_k^2  \rightarrow \int_{\partial \omega} q \phi^2 .$$
\end{lemma}
\textbf{Proof. }Let us first note that since $(\phi_k)_k$ is bounded in $H^1 (\omega)$, we already know that up to extraction of a subsequence, it converges weakly in $H^1 (\omega)$ and strongly in $L^2 (\omega )$ to a function $\phi$. Moreover, it follows from the traces theory that $((\phi_k)_{| \partial \omega })_k$ is bounded in~$W^{1/2,2} (\partial \omega)$. Thus, up to the extraction of some subsequence, it converges in $L^2 (\partial \omega )$ to the trace $\phi_{| \partial \omega }$. That is, we have for any $\lambda \in \mathbb{R}$, as $k \rightarrow +\infty$:
\begin{equation*}
\int_{\partial \omega} q \phi_k^2 \rightarrow  \int_{\partial \omega}q  \phi^2.
\end{equation*}
Therefore, it now remains to show that
\begin{equation*}
\int_\omega \frac{\partial h_k}{\partial T}  (.,0) \phi_k ^2  - \int_{\partial \omega} q  \phi_k^2 \rightarrow 0 .
\end{equation*}
From the hypothesis (\ref{eqn:hypbord}) and the $L^2$-bound on $(\phi_k)_{k \in \mathbb{N}}$, by noting $\Gamma_{\varepsilon_k}\!=\!\omega \cap (\partial \omega + \overline{B} (0,\varepsilon_k))$, we only have to prove that:
\begin{equation}\label{eqn:etape1b}
\int_{\Gamma_{\varepsilon_k}} \frac{\partial h_k}{\partial T}  (.,0) \phi_k ^2  - \int_{\partial \omega} q \phi_k^2 \rightarrow 0 . \end{equation}
Let the function $d$ be the distance from the boundary $\partial \omega$. It follows from the coarea formula that for $\varepsilon_k$ small enough (that is, for a sufficiently large $k$), we have:
$$\displaystyle \int_0^{\varepsilon_k} \displaystyle  \left( \int_{d^{-1} (s) \cap \Gamma_{\varepsilon_k} } \frac{\partial h_k}{\partial T}  (.,0) \phi_k^2  \right) ds = \int_{\Gamma_{\varepsilon_k}} \frac{\partial h_k}{\partial T}  (.,0) \phi_k ^2 .$$
For $\varepsilon_k$ small enough and $0 \leq s \leq \varepsilon_k$, we can parametrize $d^{-1} (s) \cap \Gamma_{\varepsilon_k} $ by $z-s n(z)$, where~$z \in \partial \omega$ and $n(z)$ is the outward normal unit of $\partial \omega$ on $z$. We then obtain:
\begin{eqnarray*}
& \displaystyle (1 + O (\varepsilon_k )) \int_0^{\varepsilon_k}  \left( \int_{d^{-1} (s) \cap \Gamma_{\varepsilon_k} } \frac{\partial h_k}{\partial T}  (.,0) \phi_k^2  \right) ds \vspace{3pt} \\
  = \displaystyle & \displaystyle \int_0^{\varepsilon_k}   \int_{\partial \omega} \frac{\partial h_k}{\partial T}  (z-s n(z),0)  \phi_k( z-s n(z))^2  dz ds \vspace{3pt} \\
 = \displaystyle  & \qquad \displaystyle \int_0^1   \int_{\partial \omega} \varepsilon_k \frac{\partial h_k}{\partial T}  (z-\varepsilon_k s n(z),0)  \phi_k( z- \varepsilon_k s n(z))^2  dz ds . \end{eqnarray*}
We then have on one hand:
\begin{eqnarray*}
\begin{array}{c}
\displaystyle \int_0^1   \int_{\partial \omega} \varepsilon_k \frac{\partial h_k}{\partial T}  (z-\varepsilon_k s n(z),0) \left( \phi_k( z- \varepsilon_k s n(z))^2 - \phi_k( z)^2 \right)dz ds \vspace{3pt} \\
 \displaystyle = -\int_0^1   \int_{\partial \omega} \varepsilon_k \frac{\partial h_k}{\partial T}  (z-\varepsilon_k s n(z),0)  \int_0^s 2 \varepsilon_k \phi_k (z-\varepsilon_k \tau n(z)) \nabla \phi_{k} (z-\varepsilon_k \tau n(z)).n(z) d\tau dz ds .
 \end{array}
\end{eqnarray*}
Thus, with the hypothesis (\ref{eqn:hypO}) on $h$:
\begin{eqnarray*}
\left|  \int_0^1   \int_{\partial \omega} \varepsilon_k \frac{\partial h_k}{\partial T}  (z-\varepsilon_k s n(z),0) ( \phi_k( z- \varepsilon_k s n(z))^2 - \phi_k( z)^2 )  dz ds \right|\\
\leq C \int_0^1   \int_{\partial \omega} \int_0^s 2 \varepsilon_k | \phi_k (z-\varepsilon_k \tau n(z)) \nabla \phi_{k } (z-\varepsilon_k \tau n(z)).n(z) | d\tau dz ds .
\end{eqnarray*}
From the coarea formula, with the notation $\Gamma_{\varepsilon_k ,s} = (\partial \omega + \overline{B}(0,\varepsilon_k s)) \cap \omega$, we then obtain:
\begin{equation}\label{eqn:cv1}
\begin{array}{l}
\displaystyle  \left|  \int_0^1   \int_{\partial \omega} \varepsilon_k \frac{\partial h _k}{\partial T}  (z-\varepsilon_k s n(z),0)  ( \phi_k( z- \varepsilon_k s n(z))^2 - \phi_k( z)^2 ) dz ds  \right| \vspace{3pt}\\
\displaystyle  \leq  C (1 + O (\varepsilon_k )) \int_0^1  \left( \int_{\Gamma_{\varepsilon_k, s}}  \ | \phi_k | \ \ \| \nabla \phi_{k} \| \right) ds \vspace{3pt}\\
 \displaystyle \leq  \displaystyle C' \|  \phi_k \|_{L^2 (\Gamma_{\varepsilon_k ,1})} \| \nabla \phi_k \|_{L^2 (\Gamma_{\varepsilon_k ,1})} \vspace{7pt}\\
 \displaystyle \leq \displaystyle C'' \|  \phi_k \|_{L^2 (\Gamma_{\varepsilon_k , 1})} \vspace{3pt}\\
 \rightarrow \displaystyle 0 .
\end{array}
\end{equation}
Here, we used the fact that the sequence $(\phi_k)_k$ is bounded in $H^1 (\omega)$, and converges strongly in $L^2 (\omega )$. On the other hand, it immediately follows from (\ref{eqn:hyp}) that
\begin{eqnarray}\label{eqn:cv2}
\int_0^1   \int_{\partial \omega} \varepsilon_k \frac{\partial h_k}{\partial T}  (.,0) (z-\varepsilon_k s n(z)) \phi_k( z)^2  dz ds - \int_{\partial \omega} q \phi_k^2 \rightarrow 0 .
\end{eqnarray}
Then, (\ref{eqn:cv1}) and (\ref{eqn:cv2}) imply (\ref{eqn:etape1b}), which concludes the proof of Lemma \ref{lemma1}.
$\Box$

\subsection{Locally uniform convergence of $\mu_{h_k}$ to $\nu_q $}\label{sec:cvmunu}

We now begin the proof of Theorem \ref{th:main}. Let us fix $\lambda \in \mathbb{R}$. It follows from (\ref{eqn:principaleigenvalueinfmu}) that
$$\mu_{h_k} (\lambda ) \leq \int_\omega | \nabla \psi_{q,\lambda} (y) |^2 dy - \lambda \int_\omega u(y) \psi_{q,\lambda} ^2 (y) dy + \int_\omega \left( \frac{\partial h_k}{\partial T}  (y,0)  -\frac{\partial f}{\partial T} (y,0)\right) \psi_{q,\lambda} ^2 (y) dy,$$
where $\psi_{q, \lambda}$ is the principal eigenfunction of (\ref{eqn:principaleigenvaluerobin}) normalized so that $\| \psi_{q, \lambda} \|_{L^2 (\omega )} =1$. Moreover, by Lemma~\ref{lemma1} and since $\psi_{q,\lambda}~\in~H^1~(\omega)$, we have that
$$\int_\omega \frac{\partial h_k}{\partial T}  (.,0)  \psi_{q,\lambda} ^2  \rightarrow \int_{\partial \omega} q  \psi_{q,\lambda} ^2 .$$
Thus, by passing to the limit and using (\ref{eqn:principaleigenvalueinfnu} ), we obtain that for all $\lambda \in \mathbb{R}$:
\begin{eqnarray}\label{eqn:munbounded1}\displaystyle
& & \displaystyle \limsup \mu_{h_k} (\lambda ) \vspace{3pt} \nonumber \\
 & & \leq \displaystyle \int_\omega | \nabla \psi_{q,\lambda} (y) |^2 dy - \lambda \int_\omega u(y) \psi_{q,\lambda} ^2 (y) dy + \int_{\partial \omega} q \psi_{q,\lambda}^2 - \int_\omega \frac{\partial f}{\partial T} (y,0) \psi_{q,\lambda}^2 (y) dy \vspace{3pt}  \nonumber \\
\displaystyle
 & & = \displaystyle \nu_q  (\lambda) .
\end{eqnarray}
We can also deduce that the sequence $(\mu_{h_k} (\lambda))_k$ is bounded for all $\lambda \in \mathbb{R}$. Indeed, it is bounded from above because of (\ref{eqn:munbounded1}), and it is also bounded from below thanks to the first part of (\ref{eqn:principaleigenvalueinfmu}) (from $\frac{\partial h_k}{\partial T}  (.,0) \geq 0$ for all $k$ and since $u$ and $\frac{\partial f}{\partial T}  (.,0)$ are bounded). Thus, up to the extraction of a subsequence, we can assume that $\mu_{h_k} (\lambda)$ converges to some limit~$\mu  (\lambda)$. We now want to show that $\mu  (\lambda) = \nu_q  (\lambda)$. From (\ref{eqn:principaleigenvalueinfmu}), we have that:
\begin{eqnarray*}
\begin{array}{c}
\displaystyle
\int_\omega | \nabla \phi_{h_k,\lambda} (y) |^2 dy - \lambda \int_\omega u(y) \phi_{h_k,\lambda} ^2 (y) dy + \int_\omega \left( \frac{\partial h_k}{\partial T}  (y,0) -\frac{\partial f}{\partial T} (y,0)\right) \phi_{h_k,\lambda} ^2 (y) dy \vspace{3pt} \\ \displaystyle = \ \int_\omega  \mu_{h_k} (\lambda ) \phi_{h_k,\lambda} ^2 (y) dy.
\end{array}
\end{eqnarray*}
We know that $u$ and $\frac{\partial f}{\partial T} (.,0)$ are in $L^\infty (\omega)$. Since $\frac{\partial h_k}{\partial T}  (.,0) \geq 0$, $\| \phi_{h_k,\lambda} \|_2 =1$ and since the sequence~$(\mu_{h_k} (\lambda))_k$ is bounded, it then follows that for all $\lambda \in \mathbb{R}$:
$$\sup_{k \in \mathbb{N}} \ \left( \int_\omega | \nabla \phi_{h_k,\lambda} (y) |^2 dy \right) < + \infty .$$
Therefore, for all $\lambda \in \mathbb{R}$, the sequence $(\phi_{h_k,\lambda})_k$ is bounded in $H^1 (\omega )$. Up to the extraction of some subsequence, we can then assume that there exists $\phi \in H^1 (\omega )$ such that:
$$\phi_{h_k,\lambda} \rightarrow \phi \mbox{ weakly in } H^1 (\omega ) \mbox{, strongly in } L^2 (\omega).$$
We consider each term in (\ref{eqn:principaleigenvalueinfmu}) in order to pass to the limit in $k \rightarrow +\infty$:
$$ \liminf \int_\omega | \nabla \phi_{h_k,\lambda} (y) |^2 dy \geq \int_\omega | \nabla \phi(y) |^2 dy ,$$
$$- \lambda \int_\omega u(y) \phi_{h_k,\lambda} ^2 (y) dy \rightarrow - \lambda \int_\omega u(y) \phi^2 (y) dy ,$$
$$-\int_\omega \frac{\partial f}{\partial T} (y,0)) \phi_{h_k,\lambda} ^2 (y) dy \rightarrow -\int_\omega \frac{\partial f}{\partial T} (y,0)) \phi ^2 (y) dy .$$
Here, we only used the weak convergence of the sequence $(\phi_{h_k,\lambda})_k$ in $H^1 (\omega )$ and its strong convergence in $L^2 (\omega)$. Lastly, from Lemma \ref{lemma1} and up to extraction of some subsequence, $$\int_\omega \frac{\partial h_k}{\partial T}  (.,0) \phi_{h_k,\lambda} ^2 \rightarrow \int_{\partial \omega} q \phi^2 .$$
Therefore, by passing to the limit in (\ref{eqn:principaleigenvalueinfmu}):
$$
\mu  (\lambda ) \geq \int_\omega | \nabla \phi (y) |^2 dy - \lambda \int_\omega u(y) \phi ^2 (y) dy - \int_\omega \frac{\partial f}{\partial T} (y,0) \phi ^2 (y) dy + \int_{\partial \omega} q \phi ^2.
$$
From (\ref{eqn:principaleigenvalueinfnu}), it implies that $\nu_q  (\lambda ) \leq \mu  (\lambda )$ (since $\| \phi \|_2 =1$) and it is in fact an equality from (\ref{eqn:munbounded1}). By uniqueness of the limit, we have proven the simple convergence of $\mu_{h_k} (\lambda)$ toward $\nu_q  (\lambda)$.

Here, we have also shown that
$$\nu_q  (\lambda ) = \int_\omega | \nabla \phi (y) |^2 dy - \lambda \int_\omega u(y) \phi ^2 (y) dy + \int_{\partial \omega} q \phi ^2 - \int_\omega \frac{\partial f}{\partial T} (y,0)) \phi ^2 (y) dy$$
where $\phi$ is the limit, up to extraction of some subsequence, strongly in $L^2 (\omega)$ and weakly in $H^1 (\omega)$, of the sequence $(\phi_{h_k,\lambda})_k$ of the $L^2$-normalized principal eigenfunctions of pro\-blem~(\ref{eqn:principaleigenvalueh}). Thus, by nonnegativity of $\phi$ ($\phi_{h_k ,\lambda}$ nonnegative for all $k \in \mathbb{N}$, and $\| \phi\|_2 =1$) and by uniqueness of the limit, the whole sequence $(\phi_{h_k,\lambda})_k$ converges strongly in $L^2 (\omega)$ and weakly in $H^1 (\omega)$ to the $L^2$-normalized principal eigenfunction $\psi_{q , \lambda }$ of problem~(\ref{eqn:principaleigenvaluerobin}).
\\

Moreover, we remind that for all $\lambda \in \mathbb{R}$ (see (\ref{eqn:derivative})): $\mu '_{h_k} (\lambda ) = - \int_\omega u(y) \phi_{h_k,\lambda}^2 (y) dy.$
Thus~$\| \mu  '_{h_k} \|_\infty \leq \| u \|_\infty \| \phi_{h_k,\lambda} \|^2_2 = \| u \|_\infty$. It then follows from the Dini theorem that $\mu_{h_k} (\lambda)\rightarrow \nu_q  (\lambda)$ as $k\rightarrow +\infty $ uniformly on any compact subset of $\mathbb{R}$.

\begin{remark}\label{remarkphi}
The results above will be used to control the variations of the $(\phi_{h_k,\lambda})_k$ in dimension 2 in order to prove Corollary \ref{th:main22}. In fact, we will use a little more general result, where $\lambda$ is replaced by a converging sequence $(\lambda_k)_k$. This sequence would then be bounded and the $H^1$-estimates above on the eigenfunctions would still hold. One could easily check that we would then obtain the convergence of the $L^2$-normalized principal eigenfunctions of problem $(\ref{eqn:principaleigenvalueh})$ with $\lambda=\lambda_k$ toward the $L^2$-normalized principal eigenfunction of pro\-blem~$(\ref{eqn:principaleigenvaluerobin})$ with $\lambda = \lim \lambda_k$.
\end{remark}

\subsection{Convergence of the minimal speeds}

We first show the following lemma, which will be used several times throughout this paper:
\begin{lemma}\label{lemmath1}
Under the hypotheses of Theorem $\ref{th:main}$, let $c \in \mathbb{R}$, $(c_k)_{k \in \mathbb{N}}$ and $(\lambda_k)_{k \in \mathbb{N}} $ such that $c_k \rightarrow c$ and for all $k \in \mathbb{N}$, $\mu_{h_k} (\lambda_k ) = \lambda_k^2 - c_k \lambda_k .$
Then the sequence $(\lambda_k)_{k \in \mathbb{N}}$ is bounded and $\nu_q (\lambda_\infty ) = \lambda_\infty^2 -c\lambda_\infty$ for any accumulation point~$\lambda_\infty$.
\end{lemma}
\textbf{Proof of Lemma \ref{lemmath1}.} Let $(c_k)_{k \in \mathbb{N}}$ and $(\lambda_k)_{k \in \mathbb{N}} $ such that $c_k \rightarrow c \in \mathbb{R}$ and that for any~$k \in \mathbb{N}$, $\mu_{h_k} (\lambda_k ) = \lambda_k^2 - c_k \lambda_k .$
The sequence $(\lambda_k)_{k \in \mathbb{N}} $ is bounded, since by concavity of the $\mu_{h_k} $:
$$\lambda_k ^2 - c_{k} \lambda_k \leq \mu_{h_k} (0 ) + \mu_{h_k} ' (0) \lambda_k ,$$ and from the fact that the sequences $(c_{k})_k$, $(\mu_{h_k}  ' (0))_k$ and $(\mu_{h_k} (0 ))_n$ are bounded. Let now~$\lambda_\infty$ be an accumulation point of the sequence $(\lambda_n)_n$. By the uniform convergence of~$\mu_{h_k} $ on any compact, we deduce that $\lambda_\infty ^2 - c \lambda_\infty = \nu_q  (\lambda_\infty) .\  \Box$
\\

We now get back to the proof of Theorem \ref{th:main}, and assume that $\nu_q (0) < 0$. Note that this hypothesis hadn't been used in the proofs above, which thus hold whether or not the minimal speeds are well defined. By uniform convergence of $\mu_{h,k}$ toward $\nu_q$, we can also assume up to extraction of some subsequence that $\mu_{h_k} (0) < 0$ for all $k \in \mathbb{N}$. Under those assumptions, we define the minimal speeds $c_q^*$ and $c_{h_k}^*$ for any $k\in \mathbb{N}$ as in (\ref{eqn:minspeed}). We now show that $c_{h_k}^* \rightarrow c_q^*$ to conclude the proof of Theorem \ref{th:main}.

First, let $c>c_q^*$. There exist $\varepsilon >0$ and $\lambda >0$ such that $\lambda ^2 - c \lambda \leq \nu_q   (\lambda ) - \varepsilon .$ Then, for sufficiently large $k$, we have $\lambda ^2 - c \lambda \leq \mu_{h_k} (\lambda ).$ Therefore, since $\mu_{h_k} (0) <0$ for all $k$, we have that for $k$ large enough, there exists $\lambda ' >0$ such that
$${\lambda '}^2 - c \lambda '  = \mu_{h_k} (\lambda ').$$
Hence $\limsup c_{h_k}^* \leq c$ for all $c > c_q^*$, and then
\begin{eqnarray}\label{eqn:cnbounded}
\limsup c_{h_k}^* \leq c_q^*.
\end{eqnarray}

We then deduce that the sequence $(c_{h_k}^*)_k$ is bounded. Indeed, it follows from (\ref{eqn:cnbounded}) that it is bounded from above. We also have that $ c_{h_k}^* \geq -\mu_{h_k}  '(0)$ (by concavity of $\mu_{h_k} $), and thus the sequence is bounded from below (remind that the sequence $(\mu_{h_k}  ')_k$ is bounded in~$L^\infty$). We can now assume, up to the extraction of a subsequence, that $c_{h_k}^*$ converges to some~$c \in \mathbb{R}$ verifying $c \leq c_q^*$.

Let $(\lambda_k)_k$ be a sequence of positive real numbers such that $\lambda_k ^2 - c_{h_k}^* \lambda_k = \mu_{h_k} (\lambda_k )$ for all $k \in \mathbb{N}$. By Lemma \ref{lemmath1}, this sequence is bounded and up to the extraction of some subsequence, we can assume that $\lambda_k $ converges to a $\lambda \in \mathbb{R}^+$ such that $\lambda ^2 - c \lambda = \nu_q  (\lambda) .$ Besides, since $\nu_q (0) <0$, we have that $\lambda >0$. Therefore, $c_q^* \leq c$ and $c=c_q^*$. By uniqueness of the limit, we have shown that $\lim_{k \rightarrow +\infty} c_{h_k}^* = c_q^*$ for the whole sequence, and the proof of Theorem \ref{th:main} is complete. $\Box$

\section{Convergence of some solutions}\label{sec:cvsol}

We now begin the proof of Theorem \ref{th:main2}. We recall our assumptions: $(h_k)_{k \in \mathbb{N}}$ is a sequence of functions verifying (\ref{eqn:hypbord}), (\ref{eqn:hypO}), (\ref{eqn:hyp}) and (\ref{condh}) with $h=h_k$ for all $k \in \mathbb{N}$, and such that $h_k (y,.)$ is linear for all $y \in \overline{\omega}$ and $k \in \mathbb{N}$. We let $(T_k , Y_k)$ be a sequence of solutions of problem (\ref{eqn:syshfront})-(\ref{eqn:neumann}) and (\ref{eqn:condinfty}) with $h=h_k$, $c>c_q^*$, and such that $0<T_k$ and $0<Y_k<1$.

Let also $(\lambda_k)_{k \in \mathbb{N}}$ be a sequence of positive real numbers such that $\lambda_k ^2 - c \lambda_k = \mu_{h_k} (\lambda_k )$ for any $k \in \mathbb{N}$. Lastly, we assume that there exists $0 < \Lambda_1 < \Lambda_2$ and $C_1 , \ C_2 , \ C_3 >0$ such that for all~$k \in \mathbb{N}$ and $(x,y) \in \overline{\Omega}$, $(T_k ,Y_k)$ satisfy (\ref{eqn:Tnabove}) and (\ref{eqn:Tnbelow}).

Let us first note from Lemma \ref{lemmath1} that up to extraction of some subsequence,~$\lambda_k \rightarrow \lambda_\infty$ such that
\begin{equation}\label{eqn:Texpboundlambda}
\lambda_\infty^2 - c \lambda_\infty = \nu_{q} (\lambda_\infty ).
\end{equation}
Moreover, since the real numbers $\lambda_k$ are positive and since $\nu_q (0) < 0$, we have that $ \lambda_\infty >0$. This important fact will be used several times along this section and the next one. In particular, with the hypothesis (\ref{eqn:Tnabove}), it implies that the sequence $(T_k)_{k}$ is locally bounded.
\\

We also recall the following theorem from \cite{giletti1}, giving some qualitative properties of the traveling front solutions of (\ref{eqn:syshfront})-(\ref{eqn:neumann}) and (\ref{eqn:condinfty}):
\begin{theorem}\label{th:thqualit}
Let $(c,T,Y)$ be a solution of $(\ref{eqn:syshfront})$-$(\ref{eqn:neumann})$ and $(\ref{eqn:condinfty})$ such that $0<T$ and $0<Y<1$. Then $T$ is bounded, $T(-\infty ,.)=0$, $Y(-\infty ,.)=Y_{\infty} \in (0,1)$.
\end{theorem}

\subsection{Exponential bound on $(Y_k)_{k \in \mathbb{N}}$}\label{seclemmaY}

\begin{lemma}\label{lemmaY}
Under the hypotheses of Theorem $\ref{th:main2}$, there exist $\beta>0$ and $\gamma >0$ such that for any $k\in \mathbb{N}$ :
\begin{eqnarray*}
\max (0, 1- \gamma e^{- \beta x} ) \leq Y_k < 1 .
\end{eqnarray*}
\end{lemma}
\textbf{Proof. }It has already been said that $Y_k <1$ for all $k \in \mathbb{N}$, which comes from the fact that we only consider non-trivial solutions. We introduce the following principal eigenvalue problem (\ref{eqn:principaleigenvaluelemma}), depending on a parameter $\lambda \in \mathbb{R}$:
\begin{equation}\label{eqn:principaleigenvaluelemma}
\left\{
\begin{array}{rcll}
\displaystyle  -\Delta_y \chi_\lambda - \lambda u(y) \chi_\lambda & = & \rho (\lambda ) \chi_\lambda & \mbox{ in } \omega , \vspace{3pt} \\
\displaystyle \frac{\partial \chi_\lambda }{\partial n} & = & 0 & \mbox{ on } \partial \omega .\\
\end{array}
\right.
\end{equation}
This is the same principal eigenvalue problem as (\ref{eqn:principaleigenvalueh}) and (\ref{eqn:principaleigenvaluerobin}), with $q=h=f=0$ (the purpose of its introduction is only to simplify some of our notations). In particular, we have that $\rho (\lambda)$ is concave. Furthermore, (\ref{eqn:derivative}) with $h=f=0$, together with the fact that any positive constant is an eigenfunction of (\ref{eqn:principaleigenvaluelemma}) with $\lambda =0$, imply that $\rho (0) =\rho ' (0) =0<c$. The fact that $c$ is positive follows from the first part of Theorem \ref{th:thremind}, proved in \cite{giletti1}, stating that traveling front solutions only exist with  positive speeds.

One can then choose $\beta >0$ small enough so that
\begin{equation}\label{eqn:beta}
\left\{
\begin{array}{l}
0<\beta < \inf_{k\in \mathbb{N}} \lambda_k , \vspace{3pt}\\
\rho (\beta \mbox{Le})-\beta^2 +c\beta \mbox{Le} >0 .
\end{array}
\right.
\end{equation}
Note that since each $\lambda_k$ is positive and as $\lambda_k \rightarrow \lambda_\infty >0$ up to extraction of some subsequence, we indeed have that $\inf_{k \in \mathbb{N}}~\lambda_k~>0$. Let also $\gamma > 0$ large enough so that
\begin{equation}\label{eqn:gamma}
\left\{
\begin{array}{l}
\displaystyle \gamma \times \min_{\overline{\omega}} \chi_{\beta \text{Le}} \geq 1 ,\\
\displaystyle \gamma \mbox{Le}^{-1} (\rho (\beta \mbox{Le})-\beta^2 +c\beta \mbox{Le}) \times \min_{\overline{\omega}} \chi_{\beta \text{Le}} > C_1 \max_{y \in \overline{\omega}} \left( \frac{\partial f}{\partial T} (y,0) \right) ,
\end{array}
\right.
\end{equation}
where $\chi_{\beta \text{Le}}$ is the positive eigenfunction of (\ref{eqn:principaleigenvaluelemma}) with $\lambda=\beta \mbox{Le}$, normalized in such a way that $\|\chi_{\beta \text{Le}}\|_{L^\infty (\omega )}~=~1$. Let $\underline{Y}$ be defined by
\begin{eqnarray*}
\underline{Y} (x,y) = \max(0,1-\gamma \chi_{\beta \text{Le}} (y) e^{-\beta x}).
\end{eqnarray*}
Note that $\underline{Y} =0$ for $x \leq 0$. Let us check that for any $k\in \mathbb{N}$, $\underline{Y}$ is a sub-solution for (\ref{eqn:syshfront})-(\ref{eqn:neumann}) with $T={T}_k$ and $h=h_k$. Note first that $\underline{Y}$ satisfies the Neumann boundary conditions on $\partial \Omega$. Moreover, when $\underline{Y} >0$, then $x>0$ and
\begin{eqnarray*}
&& \text{Le}^{-1} \Delta \underline{Y} + (c-u(y))\underline{Y}_{x} - f(y,{T}_k)\underline{Y} \vspace{6pt} \\
&&  \geq \gamma \text{Le}^{-1} ( \rho (\beta \text{Le}) -\beta^2 +c\beta \text{Le}) \chi_{\beta \text{Le}} (y) e^{-\beta x} \vspace{6pt}
 - \ C_1 \frac{\partial f}{\partial T} (y,0)  e^{-\lambda_k x} (1 - \gamma \chi_{\beta \text{Le}} (y) e^{-\beta x}) \vspace{3pt} \\
&& \geq \gamma \text{Le}^{-1} ( \rho(\beta \text{Le}) -\beta^2 +c\beta \text{Le}) \chi_{\beta \text{Le}} (y) e^{-\beta x} - C_1 \frac{\partial f}{\partial T} (y,0) e^{-\beta x} \vspace{3pt} \\
&& \geq 0 ,
\end{eqnarray*}since $f$ of the KPP-type, and because of (\ref{eqn:beta})-(\ref{eqn:gamma}).

Besides, we have that $\underline{Y} (-\infty,.)=0 < Y_k$ and $\underline{Y} (+\infty ,.)=1=Y_k (+\infty, .)$ for each~$k\!\in\!\mathbb{N}$. Therefore, it follows from the weak maximum principle in unbounded domains that $\underline{Y} \leq Y_k$ in~$\Omega$. This concludes the proof of Lemma~\ref{lemmaY}. $\Box$

\subsection{$H^1_{loc} (\overline{\Omega})$ estimates on $(T_k ,Y_k)$}\label{sec:estim1}

For any $k \in \mathbb{N}$, $(T_k ,Y_k)$ verifies (\ref{eqn:syshfront}) with $h=h_k$, together with Neumann boundary conditions and the conditions at infinity (\ref{eqn:condinfty}).
We first integrate the equation verified by~$Y_k$ over $(-N,N) \times \omega$ where $N\in \mathbb{R}^+$. We obtain:
\begin{eqnarray*}
& \displaystyle \int_\omega \displaystyle \Big[ \mbox{Le}^{-1} \left(Y_{k,x} (N,y) - Y_{k,x} (-N,y)\right) + (c-u(y))\left(Y_k (N,y) - Y_k (-N,y)\right) \Big]dy & \\
& \displaystyle =\int_{(-N,N) \times \omega} f(y,T_k (x,y))Y_k (x,y) dxdy .&
\end{eqnarray*}
But for each $k$, the left-hand side is bounded independently of $N$ (since $0 < Y_k < 1$ and~$Y_{k,x} /Y_k$ is bounded from the Harnack inequality for each $k$) and the function $f(T_k)Y_k$ is positive, thus its integral over $\Omega$ converges. Moreover, since for all $k$, $Y_{k,x}(\pm \infty) =0$ and~$0< Y_k < 1$, we obtain by passing to the limit $N \rightarrow +\infty $:
\begin{equation}\label{majf}
 \sup_{k \in \mathbb{N}} \int_\Omega f(y,T_k (x,y))Y_k (x,y) dxdy \leq \int_\omega | c-u(y) | dy < +\infty.
\end{equation}
For any $k \in \mathbb{N}$, we multiply by $Y_k$ the equation verified by $Y_k$ and integrate over~$(-N,N)~\times~\omega$:
\begin{eqnarray*}
& \displaystyle \int_\omega \left[\mbox{Le}^{-1} (Y_{k,x} (N,y)Y_k (N,y)\!-\!Y_{k,x} (-N,y)Y_k (-N,y))\!+\!\frac{1}{2}(c-u(y))\!(Y_k^2 (N,y)\!-\!Y_k^2(-N,y))\right]dy & \vspace{3pt}\\
& \displaystyle =\int_{(-N,N) \times \omega} f(y,T_k)Y_k^2 dxdy + \mbox{Le}^{-1} \int_{(-N,N) \times \omega} | \nabla Y_k |^2 dxdy &\vspace{3pt}\\
& \displaystyle \geq \mbox{Le}^{-1} \int_{(-N,N) \times \omega} | \nabla Y_k |^2 dxdy .
\end{eqnarray*}
The left-hand side is again bounded independently of $N\in \mathbb{R}$. Thus
$$\int_{\Omega} | \nabla Y_k |^2 dxdy < +\infty,$$
for all $k \in \mathbb{N}$. By passing to the limit as $N \rightarrow +\infty$ and using $0<Y_k<1$, we even have that:
$$ \sup_{k \in \mathbb{N}} \int_{\Omega} | \nabla Y_k |^2 dxdy \leq \frac{\mbox{Le}}{2} \int_\omega | c- u(y)|dy < +\infty.$$
That is, the sequence $(\nabla Y_k)_{k \in \mathbb{N}}$ is uniformly bounded in $L^2 (\Omega )$.
\\

We now look for $H^1_{loc} (\overline{\Omega})$ estimates on the sequence $(T_k)_{k \in \mathbb{N}}$. We first recall that the sequence $T_k$ is locally bounded, that is, for any $K$ compact subset of $\overline{\Omega}$, we have:
\begin{eqnarray}\label{eqn:Tnbounded}
\sup_k \|T_k \|_{L^\infty (K)} <+\infty .
\end{eqnarray}
Indeed, this inequality immediately follows from hypothesis (\ref{eqn:Tnabove}) and the fact that the sequence $(\lambda_k)_{k\in \mathbb{N}}$ is bounded (since it converges to $\lambda_\infty >0$ such that $\lambda_\infty ^2 - c \lambda_\infty = \nu_{q} (\lambda_\infty )$).

By integrating the equation verified by $T_k$ over $(-N , N) \times \omega$ where $N \in \mathbb{R}$, we obtain:
$$
\displaystyle \int_\omega \Big[\big(T_{k,x} (N,y) - T_{k,x} (-N,y)\big) + (c-u(y))\big(T_k (N,y) - T_k (-N,y)\big)\Big]dy $$
$$\displaystyle =\int_{(-N,N) \times \omega} h_k (y,T_k (x,y))dxdy -\int_{(-N,N) \times \omega} f(y,T_k(x,y))Y_k(x,y) dxdy .$$
Recall that for each $k$, $T_{k,x}(\pm \infty,.) = T_k(+ \infty ,.) =0$. Moreover, it has been shown in \cite{giletti1} (as reminded here in Theorem \ref{th:thqualit}) that $T_k (- \infty ,.)=0$ for any $k$. It follows, by passing to the limit $N \rightarrow +\infty$, that:
$$\displaystyle \int_{\Omega} h_k (y,T_k (x,y))dxdy = \int_{\Omega} f(y,T_k(x,y))Y_k(x,y) dxdy .$$
In particular, the left integral converges. We then obtain from (\ref{majf}):
\begin{equation}\label{majh}
\sup_{k \in \mathbb{N}} \int_{\Omega} h_k (y,T_k (x,y))dxdy < +\infty .\end{equation}
Lastly, we multiply by $T_k$ the equation verified by $T_k$, and we integrate over $(-N , N) \times \omega$:
$$\int_\omega \Big[\big(T_{k,x} (N,y)T_k (N,y) - T_{k,x} (-N,y)T_k (-N,y)\big) + \frac{1}{2}(c-u(y))\big(T_k^2 (N,y) - T_k^2(-N,y)\big)\Big]dy$$
$$=\int_{(-N,N) \times \omega} h_k (y,T_k)T_k dxdy -\int_{(-N,N) \times \omega} f(y,T_k)Y_k T_k dxdy + \int_{(-N,N) \times \omega} | \nabla T_k |^2 dxdy .$$
We then integrate over $N \in (M,M+1)$ where $M >0$:
\begin{equation}\label{eqint2}
\begin{array}{ll}
\vspace{1mm}
& \displaystyle
\int_M^{M+1} \int_\omega \Big[(T_{k,x} (N,y)T_k (N,y) - T_{k,x} (-N,y)T_k (-N,y)) \Big] dydN \\
\vspace{1mm}
& \displaystyle + \int_M^{M+1} \int_\omega \Big[(c-u(y))(T_k^2 (N,y) - T_k^2(-N,y))\Big] dydN \\
\vspace{1mm}
\displaystyle = & \displaystyle \int_M^{M+1} \int_{(-N,N) \times \omega} h_k(y,T_k)T_k dxdydN - \int_M^{M+1}  \int_{(-N,N) \times \omega} f(y,T_k)Y_k T_k dxdydN \\
\vspace{1mm}
& \displaystyle + \int_M^{M+1} \int_{(-N,N) \times \omega} | \nabla T_k |^2 dxdydN.
\end{array}
\end{equation}
For any $M \in \mathbb{R}$, the left-hand side is bounded independently of $k \in \mathbb{N}$ from Fubini theorem and the fact that the sequence $(T_k)_{k \in \mathbb{N}}$ is locally uniformly bounded from~(\ref{eqn:Tnbounded}). Moreover, from~(\ref{majf}), (\ref{eqn:Tnbounded}) and (\ref{majh}), we have that
$$\sup_{k \in \mathbb{N}} \int_M^{M+1} \int_{(-N,N) \times \omega} h_k(y,T_k)T_k dxdydN \leq \sup_{k \in \mathbb{N}} \int_{(-M-1,M+1) \times \omega} h_k(y,T_k)T_k dxdy < +\infty ,$$
$$\sup_{k \in \mathbb{N}} \int_M^{M+1} \int_{(-N,N) \times \omega} f(y,T_k) Y_k T_k dxdydN \leq \sup_{k \in \mathbb{N}} \int_{(-M-1,M+1) \times \omega} f(y,T_k)Y_k T_k dxdy < +\infty .$$
We then conclude from (\ref{eqint2}) that for all $M >0$,
$$\sup_{k \in \mathbb{N}} \int_M^{M+1} \int_{(-N,N) \times \omega} | \nabla T_k |^2 dxdydN < +\infty,$$
and thus,
$$\sup_{k \in \mathbb{N}} \int_{(-M,M) \times \omega} | \nabla T_k |^2 dxdy < +\infty .$$
That is, the sequence $(T_k)_{k \in \mathbb{N}}$ is bounded in $H^1_{loc} (\overline {\Omega})$.

\subsection{Convergence toward a solution of (\ref{eqn:sysrobinfront})-(\ref{eqn:robin})}

By the estimates proved above, we can now assume, up to extraction of some subsequence, that the sequence $(T_k ,Y_k)_{k \in \mathbb{N}}$ converges to a pair of functions $(T,Y)$ weakly in $H^1_{loc} (\overline {\Omega})$ and strongly in $L^2_{loc} (\overline {\Omega})$. We now want to prove that $(T,Y)$ is a solution of the problem (\ref{eqn:sysrobinfront})-(\ref{eqn:robin}) and (\ref{eqn:condinfty}), and we will then show that it verifies the wanted properties.

Recall that for any $n$, $Y_k$ satisfies:
$$
\displaystyle \mbox{Le}^{-1} \Delta Y_k + (c-u(y)) Y_{k,x} -f(y,T_k)Y_k  =  0 \mbox{  in } \Omega ,
$$
with the Neumann boundary conditions on $\partial \Omega$. Since $T_k$ and $Y_k$ are at least locally bounded independently of $n$ (recall (\ref{eqn:Tnbounded})), since $f$ locally Lipschitz-continuous and from the convergence toward $(T,Y)$, it is straightforward to check that $Y$ is a weak solution of
$$
\displaystyle \mbox{Le}^{-1} \Delta Y + (c-u(y)) Y_{x} -f(y,T)Y  = 0 \mbox{  in } \Omega
$$
with the Neumann boundary conditions on $\partial \Omega$. Recall now that $T_k$ verifies
$$
\displaystyle \Delta T_k + (c-u(y)) T_{k,x} + f(y,T_k)Y_k - h_k (y,T_k)  =  0 \mbox{ in } \Omega
$$
with Neumann boundary conditions. Here, the parameters of the equation depend on $k$. Thus, the convergence is not straightforward, although it is true for the weak formulation, by the same method we used in the previous sections. Let $\phi \in C_c^\infty (\overline{\Omega} )$. By multiplying the above equation by $\phi$ and integrating over $\Omega$, we obtain:
$$ - \int_\Omega  \nabla T_k . \nabla \phi \  + \int_\Omega (c-u(y)) T_{k,x} \ \phi \  + \int_\Omega f(.,T_k) Y_k  \phi \  - \int_\Omega h_k (.,T_k) \phi \ =0 .$$
Since $T_k$ converges weakly in $H^1_{loc} (\Omega)$ and strongly in $L^2_{loc} (\Omega)$ to $T$, we have that
$$\int_\Omega  \nabla T_k . \nabla \phi \rightarrow \int_\Omega  \nabla T . \nabla \phi \mbox{ , and \ } \int_\Omega (c-u(y)) T_{k,x} \ \phi \rightarrow \int_\Omega (c-u(y)) T_{x} \ \phi ,$$
as $k \rightarrow +\infty$. Since the sequences $(T_k)_k$ and $(Y_k)_k$ are bounded in $L^{\infty}_{loc} (\Omega)$ and converges in~$L^2_{loc} (\Omega)$, and since $f$ is locally Lipschitz-continuous, we have that
$$\int_\Omega f(.,T_k) Y_k  \phi \rightarrow \int_\Omega f(.,T) Y  \phi ,$$
as $k \rightarrow +\infty$.
Lastly, since the functions $h_k (y,.)$ are assumed to be linear for any $y \in \overline{\omega}$, we have that
$$\int_\Omega h_k (.,T_k)  \phi = \int_\Omega \frac{\partial h_k}{\partial T} (.,0) T_k  \phi \rightarrow \int_{\partial \Omega } q T  \phi .$$
This result is similar to the Lemma \ref{lemma1} in Section \ref{sec:cveigen}, with $\omega$ replaced by $\Omega$. In fact, since~$\phi$ is compactly supported, one can easily check that the proof of Lemma \ref{lemma1} still holds in this case.
Therefore, we have that $T$ is a weak solution of (\ref{eqn:sysrobinfront}) with Robin boundary conditions. We conclude by standard estimates that $(T, Y)$ is a strong solution of the problem (\ref{eqn:sysrobinfront})-(\ref{eqn:robin}).

\subsection{Non-triviality and conditions at infinity}

It now only remains to be shown that $0< T$, $0<Y<1$ and that $T$, $Y$ verify the right conditions at infinity. Note first that $0 \leq T$ and $0\leq Y\leq 1$ from the convergence of~$(T_k,Y_k)$ toward $(T,Y)$. Moreover, it immediately follows from hypothesis (\ref{eqn:Tnbelow}) and Lemma~\ref{lemmaY} that there exists $(x,y) \in \Omega$ such that $T(x,y) >0$ and $Y(x,y)>0$. Thus, by the strong maximum principle, we have that $T>0$ and $Y>0$ everywhere. We now assume that~$Y=1$ somewhere. Again by the strong maximum principle, we then have that $Y=1$ everywhere. Since $Y$ is a solution of (\ref{eqn:sysrobinfront})-(\ref{eqn:robin}), it implies that $f(T)Y = 0$, and thus~$T=0$, which is a contradiction.

Let us now show that $T$ and $Y$ verify the conditions at infinity (\ref{eqn:condinfty}). It is immediate that $T(+\infty ,.) = 0$ and $Y(+\infty ,.)=1$ from the exponential bounds in (\ref{eqn:Tnabove}) and Lemma~\ref{lemmaY} (recall that the sequence $\lambda_k$ converges to some $\lambda_\infty >0$ from Lemma~\ref{lemmath1}).

In order to deal with the behavior of $(T,Y)$ on the left, we will use the next lemma:

\begin{lemma}\label{lemma2}
$\|T\|_{L^\infty (\Omega )} <+\infty$.
\end{lemma}
The proof of this lemma (which echoes a proof of \cite{hamel-nonadiabatic}) is postponed to the next subsection.
\\

By the same method as in Section \ref{sec:estim1}, one can check that
$$\int_\Omega f(y,T (x,y)) Y(x,y) dxdy < +\infty ,$$
And the integral
\begin{eqnarray}\label{eqn:intdivY}
\int_{\Omega} |\nabla Y(x,y)|^2 dxdy < +\infty
\end{eqnarray}
converges. Let now $(x_j)_{j \in \mathbb{N}}$ be any sequence such that $x_j \rightarrow -\infty$ as $j \rightarrow +\infty$. We define the functions $Y_j (x,y) = Y(x+x_j ,y)$ for each $j \in \mathbb{N}$. It follows from standard elliptic estimates and the fact that $\|T\|_{L^{\infty} (\Omega )} < +\infty$ that this sequence is bounded in $W_{loc}^{2,p} (\overline{\Omega} )$ for all $1 \leq p < +\infty$. Therefore, up to extraction of a subsequence, it converges in $C_{loc}^{1} (\overline{\Omega} )$ to a function $Y_\infty$. Because of (\ref{eqn:intdivY}), we know that $Y_\infty$ is a constant. Hence, $Y_x (-\infty ,.) =0$.

Similarly, we now integrate equation (\ref{eqn:sysrobinfront}) verified by $T$ over $(-N ,N) \times \omega$ where $N >0$, and we obtain:
\begin{eqnarray*}
& \displaystyle \int_\omega \Big[\big(T_{x} (N,y) - T_{x} (-N,y)\big) + (c-u(y))\big(T (N,y) - T (-N,y)\big)\Big]dy & \\
& \displaystyle = \int_{(-N,N) \times \partial \omega } qT(x,y)dxdy - \int_{(-N,N) \times \omega} f(y,T (x,y))Y (x,y) dxdy .
\end{eqnarray*}
By Lemma \ref{lemma2} and the Harnack inequality, we know that the left-hand side of this equation is bounded independently of $N$, and that, by passing to the limit as $N \rightarrow +\infty$:
$$\int_{\partial \Omega } qT(x,y)dxdy < +\infty .$$
Furthermore, by multiplying the equation (\ref{eqn:sysrobinfront}) satisfied by $T$ by $T$ itself, and integrating over the domain $(-N,N) \times \omega$ with $N>0$, we obtain:
\begin{equation*}
\begin{array}{c} \displaystyle \int_\omega \Big[ \big(T_x (N,y) T(N,y) - T_x (-N,y)T(-N,y)\big) + \frac{1}{2} (c-u(y))\big(T^2 (N,y) -T^2 (-N,y)\big)\Big]dy \\
\\
=\displaystyle \int_{(-N,N) \times \partial \omega} q T(x,y)^2 dxdy + \int_{(-N,N) \times \omega} \Big[|\nabla T | ^2 - f(y,T(x,y))Y (x,y) T (x,y) \Big]dx dy.
\end{array}
\end{equation*}
We conclude that the integral
\begin{eqnarray*}
\int_{\Omega} |\nabla T(x,y)|^2 dxdy < +\infty
\end{eqnarray*}
converges. As before, from standard elliptic estimates and since $\|T\|_{L^{\infty} (\Omega )}~<~\!~+\infty$, we have that $T$ converges to a constant $T_\infty$ near $-\infty$. Hence~$T_x (-\infty ,.) =0$.

As a conclusion, $(T,Y)$ is a solution of (\ref{eqn:sysrobinfront})-(\ref{eqn:robin}) and verifies (\ref{eqn:condinfty}), which ends the proof of Theorem \ref{th:main2}. $\Box$

\subsection{Proof of Lemma \ref{lemma2}}\label{sec:lemma2}

Assume for the sake of contradiction that $T$ is not in $L^\infty (\Omega)$. Let us first note that from hypothesis (\ref{eqn:Tnabove}) and Lemma \ref{lemmath1}, we know that:
\begin{eqnarray}\label{eqn:Texpbound}
0 \leq T(x,y) \leq C_1 e^{-\lambda_\infty x}
\end{eqnarray}
for all $(x,y)\in \Omega$, and where $\lambda_\infty >0$ satisfies (\ref{eqn:Texpboundlambda}). Hence, the only possibility for the function $T$ to grow is on the left, and there exists a sequence $(x_j ,y_j)_{j \in \mathbb N}$ in $\mathbb{R} \times \overline{\omega}$ so that
\begin{eqnarray}\label{eqn:lemmaabs}
T(x_j ,y_j) \rightarrow +\infty \mbox{ and } x_j \rightarrow -\infty \mbox{ as } j\rightarrow +\infty.
\end{eqnarray}
We now want to show that $Y(-\infty ,.) =0$. Since the function $| \nabla T | / T$ is globally bounded from standard elliptic estimates and the Harnack inequality up to the boundary, it follows that for each $R >0$,
$$\min_{(x,y) \in [x_j -R ,x_j +R] \times \overline{\omega}} T (x ,y) \rightarrow +\infty$$
as $j \rightarrow + \infty$. Let also $m=\min_{y \in \overline{\omega}} f(y,1) >0$. We use again the principal eigenvalue problem (\ref{eqn:principaleigenvaluelemma}), which we introduced in Section~\ref{seclemmaY}. As mentioned before, the function $\rho$ is concave and $\rho (0)  =0$. Therefore, there exist exactly two real numbers $\alpha_\pm$ such that $\alpha_{-} < 0 < \alpha_{+}$ and
$$\text{Le}^{-1} \rho (-\alpha_{\pm} \text{Le}) = \text{Le}^{-1} \alpha_{\pm}^2 + c\alpha_{\pm} -m .$$
We denote by $\chi_{\pm}$ the two principal eigenfunctions of problem (\ref{eqn:principaleigenvaluelemma}) with $\lambda~=~- \alpha_{\pm} \text{Le}$, normalized so that $\min_{\overline{\omega}} \chi_{\pm} = 1$. The functions $u_{\pm} (x,y) = e^{\alpha_{\pm} x} \chi_{\pm} (y)$ then satisfy
$$
\left\{
\begin{array}{rcll}
\displaystyle \text{Le}^{-1} \Delta u_{\pm} + (c-u(y)) u_{\pm ,x} - m u_{\pm} & = & 0 & \mbox{in } \Omega ,\vspace{3pt} \\
\displaystyle \frac{\partial u_{\pm}}{\partial n} & = & 0 & \mbox{on } \partial \Omega .\\
\end{array}
\right.
$$
Fix now any $R >0$ and choose $N \in \mathbb{N}$ so that
$$\min_{(x,y) \in [x_j -R ,x_j +R] \times \overline{\omega}} T (x ,y) \geq 1$$
for all $j \geq N$. Then, as the function $f(y,T)$ is increasing in the variable $T$, we have that~$f(y,T) \geq f(y,1) \geq m$ in $[x_j -R ,x_j +R] \times \overline{\omega}$ for all $y \in \omega$ and $j \geq N$. Hence, on the same domain,
$$\text{Le}^{-1} \Delta Y + (c-u(y)) Y_x - mY \geq 0 .$$
The function $Y$ also satisfies the Neumann boundary conditions on $\partial \Omega $. Furthermore, $Y\leq1$ in $\Omega$. It then follows from the weak maximum principle that
$$\forall (x,y) \in [x_j -R ,x_j +R] \times \overline{\omega} \mbox{,   } \ \ Y(x,y) \leq e^{\alpha_+ (x-x_j -R)} \chi_+ (y) + e^{\alpha_- (x-x_j +R)} \chi_- (y).$$
Therefore, along the section $x = x_j$, the function $Y$ is small:
$$\limsup_{j\rightarrow +\infty} \left( \max_{y \in \overline{\omega}} Y (x_j ,y) \right) \leq \max \left( \max_{\overline{\omega}} \chi_+ , \max_{\overline{\omega}} \chi_- \right) \times ( e^{- \alpha_+ R} + e^{\alpha_- R}) .$$
Since $R>0$ can be chosen arbitrary, one concludes that $Y(x_j ,.) \rightarrow 0$ uniformly in $\overline{\omega}$ as~$j \rightarrow +\infty$. Let now $\varepsilon > 0$ be any positive real number, and $N \in \mathbb{N}$ such that $Y(x_j ,y) \leq \varepsilon$ for all~$j \geq N$ and $y \in \overline{\omega}$. Since the function $Y$ satisfies
$$\mbox{Le}^{-1} \Delta Y + (c-u(y))Y_x = f(y,T) Y \geq 0 ,$$
it follows from the weak maximum principle that $Y(x,y) \leq \varepsilon$ for all $(x,y) \in [x_j ,x_N] \times \overline{\omega}$ and $j \geq N$ such that $x_j \leq x_N$. Since $x_j \rightarrow -\infty$ as $j \rightarrow +\infty$, we have that $Y \leq \varepsilon$ in~$(-\infty , x_N ] \times \overline{\omega}$. Thus, as $Y \geq 0$, $Y(-\infty,.) =0$ uniformly in $y \in \overline{\omega}$.

\par
We now use this to find an increasing exponential bound on the temperature to control its behaviour on the left, and thus to reach a contradiction with (\ref{eqn:lemmaabs}). From (\ref{eqn:Texpboundlambda}) and~(\ref{eqn:principaleigenvalueinfnu}):
\begin{equation}\label{eqn:lambda1}
 \lambda_\infty^2  -c\lambda_\infty = \nu_q (\lambda_\infty ) \leq \nu_{q,f=0} (\lambda_\infty ) - \min_{y \in \overline{\omega}} \frac{\partial f}{\partial T} (y,0) ,
\end{equation}
where $\nu_{q,f=0}$ is defined as the principal eigenvalue of (\ref{eqn:principaleigenvaluerobin}) where $f$ is replaced by zero. Let $$\varepsilon = \min \left( \frac{\nu_{q,f=0} (0)}{2} , \frac{1}{2} \min_{y \in \overline{\omega}} \frac{\partial f}{\partial T} (y,0) \right) > 0 .$$ The positivity of $\nu_{q,f=0} (0)$ is easily verified from (\ref{eqn:principaleigenvalueinfnu}) with $f=0$. Let $A\geq 0$ so that
$$\forall x \leq -A, \ \forall y \in \overline{\omega} , \ \frac{\partial f}{\partial T} (y,0)Y(x,y) \leq \varepsilon .$$
Such a $A$ exists since $Y(-\infty ,.)=0$. As a consequence of the continuity of $\nu_q$ and (\ref{eqn:lambda1}), there exists $\Lambda > \lambda_\infty$ such that
\begin{equation}\label{eqn:Lambda1}
-\nu_{q,f=0} (\Lambda ) - c\Lambda + \Lambda ^2 < - \frac{1}{2} \min_{y \in \overline{\omega}} \frac{\partial f}{\partial T} (y,0) \leq -\varepsilon .
\end{equation}
We denote by $U$ the positive function defined by
$$T(x,y)=U(x,y) e^{-\Lambda x} \psi_{f=0,\Lambda} (y),$$
where $\psi_{f=0,\Lambda}$ is the principal eigenfunction of (\ref{eqn:principaleigenvaluerobin}) with the parameter $\Lambda $ and $f=0$, normalized so that $\| \psi_{f=0,\Lambda} \|_{L^2 (\omega )} =1$. Besides, one has that $T(x,y) \leq C_1 e^{-\lambda_\infty x}$ for all~$x \leq 0$ (see (\ref{eqn:Texpbound})) and thus $U(-\infty ,.)=0$. It is also easy to verify that we have $\partial_n U =0$ on $\partial \Omega$. Furthermore, one can check that
\begin{eqnarray*}
\Delta U + (c-u(y) - 2\Lambda )U_x +2\frac{\nabla_y \psi_{f=0,\Lambda} }{\psi_{f=0,\Lambda } }.\nabla_y U & \vspace{3pt} \\
+(\Lambda^2 -\nu_{q,f=0} (\Lambda )  -c\Lambda + g(x,y))U  & = & 0 \ \ \mbox{ in } \Omega ,
\end{eqnarray*}
where
$$g(x,y) = \frac{f(y,T(x,y))}{T(x,y)} Y(x,y) \leq \frac{\partial f}{\partial T} (y,0) Y(x,y) \leq \varepsilon$$
for all $(x,y) \in (-\infty , -A]\times \omega$.  Therefore, we have
$$\Delta U +(c-u(y)-2\Lambda )U_x +2\frac{\nabla_y \psi_{f=0,\Lambda} }{\psi_{f=0,\Lambda} }.\nabla_y U + (\Lambda^2 -\nu_{q,f=0} (\Lambda )  -c\Lambda + \varepsilon )U \geq 0$$
for all $(x,y) \in (-\infty , -A]\times \omega$.

Because of (\ref{eqn:Lambda1}), we shall now apply the maximum principle to the previous operator, and look for a suitable super-solution. Since $\varepsilon \leq \nu_{q,f=0} (0) /2 < \nu_{q,f=0} (0)$, there exists $\delta >0$ such that
$$\delta^2 + c\delta - \nu_{q,f=0} (-\delta ) + \varepsilon < 0 .$$
One can then check that the function
$$\overline{U} (x,y) = e^{(\Lambda + \delta )x} \times \frac{\psi_{f=0,-\delta } (y)}{\psi_{f=0, \Lambda} (y)} \ ,$$ where $\psi_{f=0,- \delta}$ is the principal eigenfunction of (\ref{eqn:principaleigenvaluerobin}) with the parameter $-\delta $ and $f=0$,
satisfies
\begin{equation*}
\begin{array}{c}
\vspace{2mm}
\displaystyle
\Delta \overline{U} + (c-u(y)-2\Lambda ) \overline{U}_x +2 \frac{\nabla_y \psi_{f=0, \Lambda}}{\psi_{f=0, \Lambda}} . \nabla_y \overline{U} + (\Lambda ^2 - \nu_{q,f=0} (\Lambda ) -c \Lambda + \varepsilon )\overline{U} \\
\displaystyle
=(\delta^2 +c\delta -\nu_{q,f=0} (-\delta )+\varepsilon )\overline{U} \leq 0 \mbox{ in } \Omega ,
\end{array}
\end{equation*}
along with Neumann boundary conditions. It follows from the maximum principle that the difference $\overline{U} -U$ can not attain an interior negative minimum. Moreover, $\overline{U} > 0$ and one can normalize the function $\psi_{f=0,-\delta}$ so that~$U(-A,y) \leq \overline{U} (-A ,y)$ for all $y \in \overline{\omega}$. Finally, both $U$ and $\overline{U}$ tend to $0$ as $x\rightarrow -\infty$. We conclude that
$$\forall x \leq -A, \ \forall y \in \overline{\omega}, \ U(x,y)\leq \overline{U} (x,y) . $$
In other words,
$$\forall x \leq -A, \ \forall y \in \overline{\omega}, \ T(x,y)\leq e^{\delta x} \psi_{f=0,-\delta} (y) \leq \gamma e^{\delta x},$$
where $\gamma = \max_{y\in \overline{\omega}} \psi_{f=0, -\delta} (y)$, and we have reached a contradiction with (\ref{eqn:lemmaabs}). Therefore, the proof of Lemma \ref{lemma2} is complete. $\Box$

\section{Convergence of some solutions in dimension 2}\label{sec:cvsol2}

In this Section, we begin the proof of Corollary \ref{th:main22}. Our aim is to find a suitable sequence of solutions of (\ref{eqn:syshfront})-(\ref{eqn:neumann}) and (\ref{eqn:condinfty}) with $c> \max (0 ,c_q^* )$ and $h= h_k$ such that it verifies the assumptions of Theorem \ref{th:main2}. As we said in the Introduction, the construction of this sequence will echo the proof which was used in \cite{giletti1} to prove the existence of solutions of~(\ref{eqn:syshfront})-(\ref{eqn:neumann}) and (\ref{eqn:condinfty}) for $c>c_h^*$. First, we will recall the sketch of this proof. We will then show how it allows us, in dimension 2, to obtain Corollary \ref{th:main22}.

\subsection{Construction of solutions of (\ref{eqn:syshfront})-(\ref{eqn:neumann}) and (\ref{eqn:condinfty}) - \cite{giletti1}}\label{sec:remind}

We fix here $k \in \mathbb{N}$. We remind the construction of a solution of (\ref{eqn:syshfront})-(\ref{eqn:neumann}) and (\ref{eqn:condinfty}) with~$h=h_k$ and $c>c_{h_k}^*$, which is possible for $k$ large enough since $c>c_q^*$ and because of Theorem \ref{th:main}. The first step, and the only one we will detail here, is to construct sub and super-solutions of (\ref{eqn:syshfront})-(\ref{eqn:neumann}). Then, the use of a fixed point theorem on bounded cylinders allowed us in \cite{giletti1} to construct approximate solutions. Lastly, by passing to the limit in the infinite cylinder, we could obtain a solution of (\ref{eqn:syshfront})-(\ref{eqn:neumann}) with the desired qualitative properties. This is in fact a standard procedure which has also already been applied to show the existence of fronts in \cite{hamel-quenching,berestycki3,hamel-adiabatic}, which is why only the construction of sub and super-solutions will be detailed here. We refer the reader to \cite{giletti1} for the end of the proof, which will be summed up here by a lemma.

\subsubsection*{Supersolutions for Y and T}

Note first that the constant 1 is a super-solution for $Y$.

We then construct a super-solution for the $T$-equation (\ref{eqn:syshfront}) with $Y=1$.
Since $\lim_{k\rightarrow +\infty} c_{h_k}^* = c_q^* < c$, we can assume, as already underlined, that $c_{h_k}^* < c$. Hence, let~$\lambda_k$ be the smallest positive root of $\lambda^2 - \mu_{h_k} (\lambda ) =c \lambda$, and $\overline{T}_k$ be the function defined in~$\overline{\Omega}$ by
\begin{eqnarray}\label{eqn:overT}
\overline{T}_k (x,y)=\phi_{\lambda_k} (y) e^{-\lambda_k x} >0 .
\end{eqnarray}
Here $\phi_{\lambda_k}$ is the positive principal eigenfunction of (\ref{eqn:principaleigenvalueh}) with $h=h_k$ and $\lambda =\lambda_k$, normalized so that $\| \phi_{\lambda_k} \|_{L^{2} (\omega )} =1$. The function $\overline{T}_k$ satisfies the Neumann boundary conditions on~$\partial \Omega$, and is a super-solution for the equation on $T$ in (\ref{eqn:syshfront}) with $Y=1$, i.e
$$\Delta \overline{T}_k + (c-u(y)) \overline{T}_{k,x} +f(y,\overline{T}_k) -h_k (y,\overline{T}_k)$$
$$\leq \Delta \overline{T}_k + (c-u(y)) \overline{T}_{k,x} +\left(\frac{\partial f}{\partial T} (y,0)-\frac{\partial h_k}{\partial T} (y,0)\right) \overline{T}_k =0 \mbox{ in } \overline{\Omega}.$$

\subsubsection*{Sub-solution for Y}\label{sec:subsolY}

The method we use here is the same than in the proof of Lemma~\ref{lemmaY} in Section~\ref{seclemmaY}. We define $\rho$ the principal eigenvalue of (\ref{eqn:principaleigenvaluelemma}). As before, we choose $\beta >0$ which satisfies (\ref{eqn:beta}). We also let $\gamma_k > 0$ large enough so that
\begin{equation}\label{eqn:gammak}
\left\{
\begin{array}{l}
\displaystyle \gamma_k \times \min_{\overline{\omega}} \chi_{\beta_k \text{Le}} \geq 1 ,\\
\displaystyle \gamma_k \mbox{Le}^{-1} (\rho (\beta_k \mbox{Le})-\beta_k^2 +c\beta_k \mbox{Le}) \times \min_{\overline{\omega}} \chi_{\beta_k \text{Le}} > \max_{y \in \overline{\omega}} \left( \frac{\partial f}{\partial T} (y,0) \phi_{\lambda_k} (y) \right) ,
\end{array}
\right.
\end{equation}
where $\chi_{\beta_k \text{Le}}$ is the positive eigenfunction of (\ref{eqn:principaleigenvaluelemma}) with $\lambda=\beta_k \mbox{Le}$, normalized in such a way that $\|\chi_{\beta_k \text{Le}}\|_{L^\infty (\omega )}~=~1$.
\begin{remark} Unlike in Section~$\ref{seclemmaY}$, $\gamma_k$ indeed depends on $k$, since we a priori lack for estimates on $(\phi_{\lambda_k})_{k \in \mathbb{N}}$.
\end{remark}
Let $\underline{Y}_k$ be defined by
\begin{eqnarray}\label{eqn:underY}
\underline{Y}_k (x,y) = \max(0,1-\gamma_k \chi_{\beta \text{Le}} (y) e^{-\beta x}).
\end{eqnarray}
As in Section \ref{seclemmaY}, one can check that $\underline{Y}_k$ is a sub-solution for (\ref{eqn:syshfront})-(\ref{eqn:neumann}) with~$T=\overline{T}_k$ and $h=h_k$. That is, $\underline{Y}_k$ satisfies the Neumann boundary conditions on $\partial \Omega$, and for any~$(x,y)\in \Omega$:
$$ \text{Le}^{-1} \Delta \underline{Y}_k + (c-u(y))\underline{Y}_{k,x} - f(y,\overline{T}_k)\underline{Y}_k \geq 0 .$$

\subsubsection*{Sub-solution for T}\label{sec:subsolT}

Lastly, we will construct a sub-solution for $T$ with $Y=\underline{Y}_k$. Recall that $\lambda_k^2 - \mu_{h_k} (\lambda )=c\lambda_k $. We define $a_k (\lambda) = \lambda^2 - \mu_{h_k} (\lambda )$. We then show that $a_k '(\lambda_k ) < c$. Indeed, since $a_k (0)>0$ and $\lambda_k$ is the smallest positive root of $\lambda^2 - \mu_{h_k} (\lambda )=c\lambda$, we have $a_k '(\lambda_k )\leq c$. Furthermore, if $a_k '(\lambda_k) =c$, then $\lambda^2 - \mu_{h_k} (\lambda )\geq c\lambda $ for all $\lambda \in \mathbb{R}$ by convexity of $a_k$, whence $c_{h_k}^* \geq c$, which is a contradiction. We conclude, as announced, that $a_k ' (\lambda_k ) < c$.

The above allows us to choose $\eta_k>0$ small enough so that
\begin{equation}\label{eqn:eta}
\left\{
\begin{array}{l}
0<\eta_k< \min (\beta ,\alpha \lambda_k) ,\\
\varepsilon_k :=c(\lambda_k +\eta_k)-a_k(\lambda_k + \eta_k)>0 ,
\end{array}
\right.
\end{equation}
where $\alpha >0$ such that $f(y,.)$ is of class $C^{1,\alpha } ([0,s_0 ])$ for some $s_0 >0$ uniformly in $y \in \overline{\omega}$. Let $M \geq 0$ such that
\begin{equation}\label{eqn:eqM}
\begin{array}{l}
\displaystyle f(y,s)\geq \frac{\partial f}{\partial T} (y,0)s - Ms^{1+\alpha } ,\\
\end{array}
 \mbox{ for all } s \in [0,s_0]  \mbox{ and for all } y \in \overline{\omega}.
\end{equation}
Now take $x_k \geq 0$ sufficiently large so that
$$\underline{Y}_k (x,y) = 1 -\gamma_k \chi_{\beta \text{Le}} (y) e^{-\beta x} \mbox{ for all } (x,y) \in (x_k ,+\infty ) \times \overline{\omega} .$$
Next, let $\delta_k >0$ large enough so that
\begin{equation}\label{eqn:delta1}
\left\{
\begin{array}{l}
\displaystyle \phi_{\lambda_k} (y) e^{-\lambda_k x} - \delta_k \phi_{\lambda_k + \eta_k} (y) e^{-(\lambda_k +\eta_k) x} \leq s_0 \mbox{ in } \overline{\Omega} , \vspace{3pt} \\
\displaystyle \phi_{\lambda_k} (y) e^{-\lambda_k x} - \delta_k \phi_{\lambda_k + \eta_k} (y) e^{-(\lambda_k +\eta_k) x} \leq 0 \mbox{ \ in } (-\infty ,x_k ] \times \overline{\omega} ,\\
\displaystyle \delta_k \varepsilon_k \times \min_{\overline{\omega}} \phi_{\lambda_k + \eta_k} \geq \max_{y \in \overline{\omega}} \left( \gamma_k \frac{\partial f}{\partial T} (y,0)\phi_{\lambda_k} (y)+ M \phi_{\lambda_k} (y)^{1+\alpha} \right),
\end{array}
\right.
\end{equation}
where $\phi_{\lambda_k + \eta_k}$ is the positive principal eigenfunction of (\ref{eqn:principaleigenvalueh}) with $h=h_k$ and $\lambda =\lambda_k + \eta_k$, normalized so that $\| \phi_{\lambda_k + \eta_k} \|_{L^{2} (\omega )} =1$.
Lastly, we define, for all $(x,y) \in \overline{\Omega}$,
\begin{eqnarray}\label{eqn:underT}
\underline{T}_k (x,y) = \max \left( 0 , \phi_{\lambda_k} (y) e^{-\lambda_k x} -\delta_k \phi_{\lambda_k + \eta_k} (y) e^{-(\lambda_k +\eta_k)x} \right) .
\end{eqnarray}
The function $\underline{T}_k$ satisfies the Neumann boundary conditions on $\partial \Omega$. Let us now check that~$\underline{T}_k$ is a sub-solution to (\ref{eqn:syshfront}) with $Y=\underline{Y}_k$. Note first that $0 \leq \underline{T}_k \leq s_0$ in $\overline{\Omega}$. Moreover, if $\underline{T}_k (x,y) > 0$, then $x>x_k \geq 0$ whence $0\leq \underline{Y}_k(x,y)=1-\gamma_k \chi_{\beta \text{Le}} (y) e^{-\beta x}$.
Then, in that case, we have:
\begin{eqnarray*}
&& \Delta \underline{T}_k + (c-u(y))\underline{T}_{k,x} + f(y,\underline{T}_k)\underline{Y}_k - h_k (y,\underline{T}_k) \vspace{3pt} \\
&& \geq \Delta \underline{T}_k + (c-u(y))\underline{T}_{k,x} - \frac{\partial h_k}{\partial T} (y, 0)\underline{T}_k
 + \left(\frac{\partial f}{\partial T} (y,0)\underline{T}_k-M\underline{T}_k^{1+\alpha}\right)\left(1-\gamma_k \chi_{\beta \text{Le}} (y) e^{-\beta x}\right) \vspace{5pt} \\
&& \geq -\delta_k ( k(\lambda_k +\eta_k) -c(\lambda_k +\eta_k))\phi_{\lambda_k +\eta_k} (y) e^{-(\lambda_k +\eta_k)x} \vspace{3pt} \\
&& \hspace{30pt} -\frac{\partial f}{\partial T} (y,0)\gamma_k \underline{T}_k\chi_{\beta \text{Le}} (y) e^{-\beta x} -M\underline{T}_k^{1+\alpha}\\
&& \geq \delta_k \varepsilon_k \phi_{\lambda_k +\eta_k} (y) e^{-(\lambda_k +\eta_k)x} - \frac{\partial f}{\partial T} (y,0)\phi_{\lambda_k} (y) \gamma_k e^{-(\lambda_k + \beta_k )x} - M \phi_{\lambda_k} (y)^{1+\alpha} e^{-\lambda_k (1+\alpha )x}\\
&& \geq  \left(\delta_k \varepsilon_k \phi_{\lambda_k +\eta_k} (y) -\frac{\partial f}{\partial T} (y,0)\phi_{\lambda_k} (y)\gamma_k - M\phi_{\lambda_k} (y)^{1+\alpha}\right)e^{-(\lambda_k + \eta_k)x} \vspace{3pt} \\
&& \geq 0 ,
\end{eqnarray*}
because of (\ref{eqn:eta}), (\ref{eqn:eqM}), (\ref{eqn:delta1}), the fact that $h_k$ is linear and since $0~<~\phi_{\lambda_k +\eta_k}~(y)$, $0~<~\chi_{\beta \text{Le}}~(y)~\leq~1$ in~$\overline{\omega}$.

\subsubsection*{End of the construction of solutions with speed $c$ of (\ref{eqn:syshfront})-(\ref{eqn:neumann}) and (\ref{eqn:condinfty})}

We sum up the end of the proof in the following theorem.

\begin{theorem}\label{th:giletti1}
Let $\overline{T}_k$, $\underline{T}_k$ and $\underline{Y}_k$ be defined as in $(\ref{eqn:overT})$, $(\ref{eqn:underT})$ and $(\ref{eqn:underY})$ where the parameters verify the above assumptions.
Then there exists a solution $(T_k ,Y_k)$ of $(\ref{eqn:syshfront})$-$(\ref{eqn:neumann})$ with $(\ref{eqn:condinfty})$, such that $\underline{T}_k \leq T_k \leq \overline{T}_k$ and $\underline{Y}_k \leq Y_k < 1$.
\end{theorem}
As we said before, this proof relies on the use of a fixed point theorem in a truncated cylinder, and then on a passage to the limit, but we refer to \cite{giletti1} for the details.

\subsection{Proof of Corollary \ref{th:main22}}

We now assume that we are in dimension 2 ($d=2$). We want to construct a sequence~$(T_k ,Y_k)$ of solutions of (\ref{eqn:syshfront})-(\ref{eqn:neumann}) where $c> c_q^*$ and with the conditions (\ref{eqn:condinfty}), such that it verifies the assumptions (\ref{eqn:Tnabove}) and (\ref{eqn:Tnbelow}) of Theorem \ref{th:main2}.

That is, we want to find $\overline{T}_k$, $\underline{Y}_k$ and $\underline{T}_k$ sub- and super-solutions defined as in (\ref{eqn:overT}),~(\ref{eqn:underY}) and (\ref{eqn:underT}), such that there exist $0 < \Lambda_1 < \Lambda_2 $ and $C_1 , \ C_2 , \ C_3 >0$ such that for all $k \in \mathbb{N}$ and $(x,y) \in \overline{\Omega}$:
\begin{eqnarray}\label{eqn:cor1}
\overline{T}_k (x,y) \leq C_1 e^{-\lambda_k x} ,
\end{eqnarray}
\begin{eqnarray}\label{eqn:cor2}
\max (0, C_2 e^{-\Lambda_1 x} - C_3 e^{-\Lambda_2  x}) \leq \underline{T}_k (x,y).
\end{eqnarray}

We first search some $C_1 >0$ such that condition (\ref{eqn:cor1}) is satisfied. Recall (\ref{eqn:overT}):
$$\overline{T}_k (x,y)=\phi_{\lambda_k} (y) e^{-\lambda_k x} >0 ,$$
where $\phi_{\lambda_k}$ is the positive principal eigenfunction of (\ref{eqn:principaleigenvalueh}) with $h=h_k$ and $\lambda =\lambda_k$, normalized so that $\| \phi_{\lambda_k} \|_{L^{2} (\omega )} =1$. We know from Theorem \ref{th:main} that for any fixed $\lambda$, the sequence of the principal eigenfunctions of (\ref{eqn:principaleigenvalueh}) with $h=h_k$ and $L^2$-normalization is bounded in~$H^1 (\omega)$ and converges in $L^2 (\omega)$ to the principal eigenfunction $\psi_{q,\lambda}$ of (\ref{eqn:principaleigenvaluerobin}). In fact, one could easily check that this result still holds with a sequence $\lambda_k \rightarrow \lambda_\infty$ (see Remark \ref{remarkphi} in Section \ref{sec:cvmunu}). Applying this here, that means that our sequence $\phi_{\lambda_k}$ is bounded in $H^1 (\omega)$ and converges in $L^2 (\omega)$ to $\psi_{q,\lambda_\infty}$. Furthermore, since $\omega \subset \mathbb{R}$ (we are in dimension~2), we can assume, up to extraction of some subsequence, that the convergence also holds in the Holder spaces $C^{0,( 1/2) - \varepsilon} (\overline{\omega} )$ for all $\varepsilon >0$. Since $\psi_{q,\lambda_\infty}$ is a positive function (as a principal eigenfunction of (\ref{eqn:principaleigenvaluerobin})), we have that there exist $0 < K_1 < K_2$ such that for all $k \in \mathbb{N}$ large enough and $y \in \overline{\omega}$ :
\begin{eqnarray}\label{eqn:eigenfuncbound}
K_1 \leq \phi_{\lambda_k} (y) \leq K_2 .
\end{eqnarray}
With $C_1 = K_2$, the condition (\ref{eqn:cor1}) is verified.

We recall that $\beta$ is chosen as in~(\ref{eqn:beta}), and we now choose~$\gamma >0$ such that~(\ref{eqn:gammak}) holds for $\gamma =\gamma_k$, that is:
\begin{equation*}
\left\{
\begin{array}{l}
\displaystyle \gamma \times \min_{\overline{\omega}} \chi_{\beta \text{Le}} \geq 1 ,\\
\displaystyle \gamma \mbox{Le}^{-1} (\rho (\beta \mbox{Le})-\beta^2 +c\beta \mbox{Le}) \times \min_{\overline{\omega}} \chi_{\beta \text{Le}} >  K_2 \max_{y \in \overline{\omega}} \left( \frac{\partial f}{\partial T} (y,0) \right) .
\end{array}
\right.
\end{equation*}
Then:
$$\underline{Y}_k (x,y) = \underline{Y} (x,y) = \max(0,1-\gamma \chi_{\beta \text{Le}} (y) e^{-\beta x}) $$
is a suitable sub-solution for Theorem \ref{th:giletti1} for all $k$.

Lastly, we deal with condition (\ref{eqn:cor2}). We recall that we defined $a_k (\lambda) = \lambda^2 - \mu_{h_k} (\lambda )$ for any $\lambda \in \mathbb{R}$. We can also define $a (\lambda) = \lambda^2 - \nu_{q} (\lambda )$. We already know from Theorem~\ref{th:main} that $a_k \rightarrow a$ locally uniformly. Besides, as we did above for $a_k$ in Section~\ref{sec:remind}, we have that~$a (\lambda_\infty)=c\lambda_\infty$ and $a ' (\lambda_\infty ) < c$ ($\lambda_\infty$ is the smallest positive root of $a (\lambda)=c\lambda$). Let now $\eta$ small enough so that
\begin{equation*}
\left\{
\begin{array}{l}
\displaystyle 0<\eta <  \min (\beta ,\lambda_\infty ,\alpha \inf_{k \in \mathbb{N}} \lambda_k) ,\\
\displaystyle \varepsilon :=c(\lambda_\infty +\eta )-a(\lambda_\infty + \eta )>0 ,
\end{array}
\right.
\end{equation*}
where $\alpha >0$ such that $f(y,.)$ is of class $C^{1,\alpha } ([0,s_0 ])$ for some $s_0 >0$ uniformly in $y \in \overline{\omega}$. Let now $\eta_k= \lambda_\infty + \eta - \lambda_k$, which converges to $\eta$ as $k \rightarrow +\infty$. We then have for $k$ large enough that $\eta_k$ satisfies (\ref{eqn:eta}) with $\varepsilon_k \geq \frac{1}{2} \varepsilon$ bounded away from 0. Note that we used here the locally uniform convergence of $a_k$ toward $a$.
Now let $x_0 \geq 0$ sufficiently large so that
$$\underline{Y} (x,y) = 1 -\gamma \chi_{\beta \text{Le}} (y) e^{-\beta x} \mbox{ for all } (x,y) \in (x_0 ,+\infty ) \times \overline{\omega} .$$
We assume that $k$ is large enough so that:
\begin{eqnarray}\label{eqn:etapelambda}
\lambda_k \in ( \lambda_\infty- \frac{\eta}{2}, \lambda_\infty + \frac{\eta}{2}).
\end{eqnarray}
We recall that $\phi_{\lambda_k}$ converges uniformly in $\overline{\omega}$ to $\psi_{q,\lambda_\infty}$ the $L^2$-normalized positive eigenfunction of (\ref{eqn:principaleigenvaluerobin}), which implied (\ref{eqn:eigenfuncbound}). Similarly, we have that $\phi_{\lambda_k + \eta_k} =\phi_{\lambda_\infty + \eta}$ (the principal eigenfunction of (\ref{eqn:principaleigenvalueh}) with parameter $\lambda_\infty +\eta$) converges uniformly in $\overline{\omega}$ to $\psi_{q,\lambda_\infty + \eta}$ the $L^2$-normalized positive eigenfunction of (\ref{eqn:principaleigenvaluerobin}) with parameter $\lambda_\infty +\eta$. Therefore, there exist $0 < K_3 < K_4$ such that for all $k \in \mathbb{N}$ large enough and $y \in \overline{\omega}$:
\begin{eqnarray}\label{eqn:eigenfuncbound2}
K_3 \leq \phi_{\lambda_k + \eta_n} (y) \leq K_4 .
\end{eqnarray}
Since $\lambda_\infty - \eta / 2 > 0$, we can now let $\delta >0$ large enough so that
\begin{equation*}
\left\{
\begin{array}{l}
\displaystyle K_2 e^{-(\lambda_\infty + \frac{\eta}{2}) x} - \delta K_3 e^{-(\lambda_\infty +\eta ) x} \leq s_0 \mbox{ in } \overline{\Omega} , \vspace{3pt} \\
\displaystyle K_2 e^{-(\lambda_\infty - \frac{\eta}{2}) x} - \delta K_3 e^{-(\lambda_\infty +\eta ) x} \leq s_0 \mbox{ in } \overline{\Omega} ,\vspace{3pt} \\
\displaystyle K_2 e^{-(\lambda_\infty + \frac{\eta}{2}) x} - \delta K_3  e^{-(\lambda_\infty +\eta ) x} \leq 0 \mbox{ \ in } (-\infty ,x_0 ] \times \overline{\omega} ,\vspace{3pt} \\
\displaystyle K_2 e^{-(\lambda_\infty - \frac{\eta}{2}) x} - \delta K_3  e^{-(\lambda_\infty +\eta ) x} \leq 0 \mbox{ \ in } (-\infty ,x_0 ] \times \overline{\omega} ,\vspace{3pt} \\
\displaystyle \frac{ \delta \varepsilon}{2} \times K_3 \geq \gamma K_2 \max_{y \in \overline{\omega}} \left( \frac{\partial f}{\partial T} (y,0) \right) + M K_2^{1+\alpha}.
\end{array}
\right.
\end{equation*}
Thus, from (\ref{eqn:eigenfuncbound}), (\ref{eqn:etapelambda}) and (\ref{eqn:eigenfuncbound2}), $\delta$ satisfies (\ref{eqn:delta1}) with $\delta_k =\delta$ for any $k$. We can now define, for all $(x,y) \in \overline{\Omega}$ and $k \in \mathbb{N}$,
\begin{eqnarray*}
\underline{T}_k (x,y) = \max \left( 0 , \phi_{\lambda_k} (y) e^{-\lambda_k x} -\delta \phi_{\lambda_k + \eta_k} (y) e^{-(\lambda_k +\eta_k)x} \right) ,
\end{eqnarray*}
which is a suitable sub-solution for Theorem \ref{th:giletti1} for all $k \in \mathbb{N}$.

It now only remains to prove that those sub-solutions satisfy (\ref{eqn:cor2}). For $x \leq 0$, we have that $\underline{T}_k (x,y) =0$. Thus, if $\underline{T}_k > 0$, then $x > 0$ and:
\begin{eqnarray*}
\underline{T}_k (x,y) & = & \phi_{\lambda_k} (y) e^{-\lambda_k x} -\delta \phi_{\lambda_k + \eta_k} (y) e^{-(\lambda_k +\eta_k)x} \\
& \geq & K_1 e^{-(\lambda_\infty - \frac{\eta}{2}) x} -\delta K_4 e^{-(\lambda_\infty +\eta )x} \ .
\end{eqnarray*}
One can then easily conclude that
\begin{eqnarray*}
\underline{T}_k (x,y) & \geq & \max \left( 0 , K_1 e^{-(\lambda_\infty - \frac{\eta}{2}) x} -\delta K_4 e^{-(\lambda_\infty +\eta )x} \right) ,
\end{eqnarray*}
and (\ref{eqn:cor2}) is satisfied with $C_2 = K_1$, $C_3= \delta K_4$, $\Lambda_1 =\lambda_\infty - \frac{\eta}{2}$ and $\Lambda_2 = \lambda_\infty +\eta$.

Finally, our sub and super-solutions satisfy the assumptions needed for Theorem \ref{th:giletti1}, which means that there exists a sequence of solutions $(T_k ,Y_k)$ of (\ref{eqn:syshfront})-(\ref{eqn:neumann}) and (\ref{eqn:condinfty}), such that $\underline{T}_k \leq T_k \leq \overline{T}_k$ and $\underline{Y}_k \leq Y_k < 1$. Furthermore, since $\underline{T}_k$ and $\overline{T}_k$ satisfy the conditions (\ref{eqn:cor1}) and (\ref{eqn:cor2}), the sequence $(T_k ,Y_k)_{k \in \mathbb{N}}$ satisfies the assumptions (\ref{eqn:Tnabove}) and~(\ref{eqn:Tnbelow}) of Theorem \ref{th:main2}, hence the proof of Corollary~\ref{th:main22} is now complete. $\Box$

\nocite{giovangigli1}
\bibliographystyle{plain}
\bibliography{KPP5}

\begin{thebibliography}{10}

\bibitem{berestycki-hamel}
H.~Berestycki and F.~Hamel.
\newblock {\em Reaction-diffusion equations and propagation phenomena}.
\newblock Springer-Verlag, to appear.

\bibitem{hamel-quenching}
H.~Berestycki, F.~Hamel, A.~Kiselev, and L.~Ryzhik.
\newblock Quenching and propagation in {KPP} reaction-diffusion equations with
  a heat loss.
\newblock {\em Arch. Ration. Mech. Anal.}, 178:57--80, 2005.

\bibitem{berestycki5}
H.~Berestycki, F.~Hamel, and N.~Nadirashvili.
\newblock The speed of propagation for {KPP} type problems. {I} - {P}eriodic
  framework.
\newblock {\em J. European Math. Soc.}, 7:173--213, 2005.

\bibitem{berestycki6}
H.~Berestycki, F.~Hamel, and N.~Nadirashvili.
\newblock The speed of propagation for {KPP} type problems. {II} - {G}eneral
  domains.
\newblock {\em J. Amer. Math. Soc.}, 23:1--34, 2010.

\bibitem{berestycki3}
H.~Berestycki, B.~Larrouturou, and P.-L. Lions.
\newblock Multi-dimensional traveling wave solutions of a flame propagation
  model.
\newblock {\em Arch. Rational Mech. Anal.}, 111:33--49, 1990.

\bibitem{berestycki4}
H.~Berestycki and L.~Nirenberg.
\newblock Traveling wave in cylinders.
\newblock {\em Annales de l'IHP, Analyse non lin{\'e}aire}, 9:497--572, 1992.

\bibitem{giletti1}
T.~Giletti.
\newblock {KPP} reaction-diffusion equations with a non-linear loss inside a
  cylinder, preprint.

\bibitem{giovangigli1}
V.~Giovangigli.
\newblock Nonadiabatic plane laminar flames and their singular limits.
\newblock {\em SIAM J. Math. Anal.}, 21:1305--1325, 1990.

\bibitem{gordon-periodic}
P.~Gordon, L.~Ryzhik, and N.~Vladimirova.
\newblock The {KPP} system in a periodic flow with a heat loss.
\newblock {\em Nonlinearity}, 18:571--589, 2005.

\bibitem{hamel-nonadiabatic}
F.~Hamel and L.~Ryzhik.
\newblock Non-adiabatic {KPP} fronts with an arbitrary lewis number.
\newblock {\em Nonlinearity}, 18:2881--2902, 2005.

\bibitem{hamel-adiabatic}
F.~Hamel and L.~Ryzhik.
\newblock Travelling waves for the thermodiffusive system with arbitrary lewis
  numbers.
\newblock {\em Arch. Ration. Mech. Anal.}, 195:923--952, 2010.

\bibitem{kiselev1}
A.~Kiselev and L.~Ryzhik.
\newblock Enhancement of the travelling front speeds in reaction-diffusion
  equations with advection.
\newblock {\em Ann. Inst. H. Poincar{\'e} Anal. Non Lin{\'e}aire}, 18:309--358,
  2001.

\bibitem{murray1}
J.D. Murray.
\newblock {\em Mathematical biology}.
\newblock Springer, third edition, 2003.

\bibitem{roques2}
L.~Roques.
\newblock Study of the premixed flame model with heat losses the existence of
  two solutions.
\newblock {\em Euro. J. Appl. Math.}, 16:741--765, 2005.

\bibitem{williams}
F.~Williams.
\newblock {\em Combustion Theory}.
\newblock Addison-Wesley, Reading, MA, USA, 1983.

\bibitem{xin3}
J.~Xin.
\newblock Analysis and modelling of front propagation in heterogeneous media.
\newblock {\em SIAM Rev.}, 42:161--230, 2000.

\end{thebibliography}

\end{document}